\documentclass[10pt]{amsart}

\makeatletter
\@namedef{subjclassname@2020}{\textup{2020} Mathematics Subject Classification}
\makeatother

\usepackage{amsbsy, amscd, amsmath, amssymb, amstext, amsxtra, commath, eucal, float, graphicx, hyperref, latexsym, mathrsfs, mathtools, multicol, tikz, tikz-cd}

\usepackage{lineno}

\newtheorem{theorem}{Theorem}[section]
\newtheorem{lemma}[theorem]{Lemma}
\newtheorem{proposition}[theorem]{Proposition}
\newtheorem{corollary}[theorem]{Corollary}
\newtheorem{example}[theorem]{Example}
\theoremstyle{definition}
\newtheorem{definition}[theorem]{Definition}
\newtheorem{remark}[theorem]{Remark}
\newtheorem{question}[theorem]{Question}

\numberwithin{equation}{section}

\newcommand\restr[2]{\ensuremath{#1\!\!\upharpoonright_{#2}}\,}

\def\C{\mathop{\mathsf{C}\hspace{0mm}}\nolimits}
\def\cf{\mathop{\rm cf}\nolimits}

\def\CN{\mathop{\mathsf{CN}}\nolimits}

\def\int{\mathop{\rm int}\nolimits}

\def\P{\mathop{\mathsf{P}\hspace{0mm}}\nolimits}

\def\UC{\mathop{\mathsf{UC}}\nolimits}
\def\WP{\mathop{\mathsf{WP}}\nolimits}

\def\F{\mathop{\mathscr{F}\hspace{0mm}}\nolimits}

\let\mathcal\mathscr

\tikzcdset{scale cd/.style={every label/.append style={scale=#1},
    cells={nodes={scale=#1}}}}
    
\begin{document}


\title{Weakly separated spaces and Pixley-Roy hyperspaces}

\author{Alejandro R\'ios-Herrej\'on}
\address {A. R\'ios-Herrej\'on\\
Departamento de Matem\'aticas, Facultad de Ciencias, Universidad Nacional Aut\'onoma de M\'exico, Circuito ext. s/n, Ciudad Universitaria, C.P. 04510,  M\'exico, CDMX}
\email{chanchito@ciencias.unam.mx}
\thanks{The author was supported by CONACYT grant no. 814282. On behalf of all authors, the corresponding author states that there is no conflict of interest.}



\begin{abstract}
 In this paper we obtain new results regarding the chain conditions in the Pixley-Roy hyperspaces $\F[X]$. For example, if $c(X)$ and $R(X)$ denote the cellularity and weak separation number of $X$ (see Section~\ref{sec_weak_sep}) and we define the cardinals $$c^* (X) := \sup \{c(X^{n}) : n\in \mathbb{N}\} \quad \text{and} \quad R^{*}(X) := \sup \{R(X^{n}) : n\in \mathbb{N}\},$$ then we show that $R^{*}(X) = c^ {*}\left(\F[X]\right)$.
  
 On the other hand, in \cite[Question~3.23, p.~3087]{sakai2012} Sakai asked whether the fact that $\F[X]$ is weakly Lindel\"of implies that $X$ is hereditarily separable and proved that if $X$ is countably tight then the previous question has an affirmative answer. We shall expand Sakai's result by proving that if $\F[X]$ is weakly Lindel\"of and $X$ satisfies any of the following conditions: \begin{itemize}
 \item $X$ is a Hausdorff $k$-space;
 \item $X$ is a countably tight $T_1$-space;
 \item $X$ is weakly separated,
 \end{itemize} then $X$ is hereditarily separable.
\end{abstract}

\maketitle

\section{Introduction}

Chain conditions on topological spaces have been studied after their introduction by \v{S}anin in \cite{sanin1948}. Since then, many research articles have been published on the behavior of calibers, precalibers, and weak precalibers in topological spaces. For example, in \cite{rios2022}, \cite{sanin1948} and \cite{shelah1977} some results can be found regarding the preservation of these notions in the realm  of topological products.

On the other hand, the Pixley-Roy hyperspaces $\F[X]$ have also been the focus of plentiful research since their presentation by Pixley and Roy in \cite{pixroy1969}. Regarding these spaces, numerous papers and surveys have also been written (see, for example, \cite{vanDouwen1977} and \cite{lutzer1978}).

It was not until Sakai's article \cite{sakai2012} that the behavior of the {\lq\lq}precaliber $\omega_1${\rq\rq} notion in Pixley-Roy hyperspaces was studied for the first time. In that article some very interesting connections were made regarding weakly separated spaces and other cardinal functions of the hyperspace $\F[X]$.

In the present work we will delve into the study of cellularity, calibers, precalibers and weak precalibers in the hyperspaces $\F[X]$ and their relationship with other various topological concepts ($\omega$-covers, the weak separation number, the tightness, discrete families, $k$-spaces, etc).
 
\section{Preliminaries}\label{sec_prelim}

All topological and set-theoretic notions that are not explicitly mentioned in this paper should be understood as in
\cite{engelking1989} and \cite{kunen1980}, respectively. Throughout the text and unless explicitly stated otherwise, all spaces and cardinal numbers considered will be infinite.

The symbol $\omega$ will stand for both, the set of all non-negative integers and the first infinite cardinal. Additionally, the symbol $\mathbb{N}$ will stand for the set $\omega\setminus\{0\}$. If $\kappa$ is a cardinal number, the {\it cofinality} of $\kappa$ will be denoted by $\cf(\kappa)$. On the other hand, if $X$ is a set, the symbols $[X]^{<\omega}$ and $[X]^{\kappa}$ will represent the families $\{Y\subseteq X : |Y|<\omega\}$ and $\{Y\subseteq X : |Y|=\kappa\}$, respectively. Furthermore, we will denote by $\CN$ the proper class formed by the cardinal numbers, and $\UC$ will stand for the subclass of $\CN$ made up by the cardinals with uncountable cofinality.

For a topological space $X$, we will use the symbol $\tau_X$ to refer to the family of open subsets of $X$. Similarly, $\tau_X^{+}$ shall be used to denote the set $\tau_X \setminus \{\emptyset\}$ and $\tau_X^{*}$ will stand for the set $\{A\subseteq X : X\setminus A \in \tau_X\}$. Now, if $\mathcal{U}$ is a pairwise disjoint collection of $\tau_{X}^{+}$, we will say that $\mathcal{U}$ is a \textit{cellular family} in $X$.
The term {\it centered family} shall be used to designate a collection of subsets such that the intersection of the elements of any of its non-empty finite subcollections is non-empty. On the other hand, we will say that a subset $\mathcal{U}$ of $\tau_X^+$ is {\it linked} if we have $ U\cap V\neq \emptyset$ for any $U,V\in \mathcal{U}$.

A cardinal number $\kappa$ is a {\it caliber} (resp., {\it precaliber}; {\it weak precaliber}) for a topological space $X$ if for every family $\{U_\alpha : \alpha < \kappa\}\subseteq \tau_X^+$ there exists $J\in [\kappa]^{\kappa}$ such that $\{U_\alpha : \alpha \in J \}$ has non-empty intersection (resp., is centered; is linked).

We will be working throughout this paper with the following collections of cardinal numbers: \begin{align*} \WP(X)&:=\{\kappa\in\CN : \kappa \ \text{is a weak precaliber for}\ X\}; \\
\P(X)&:=\{\kappa\in\CN : \kappa \ \text{is a precaliber for} \ X\}; \ \text{and} \\
\C(X)&:=\{\kappa\in\CN : \kappa \ \text{is a caliber for} \ X\}.
\end{align*}

In the remainder of this portion of the text we will mention some auxiliary and basic propositions that we are going to use multiple times throughout the paper. Let $X$ be a topological space and $\kappa$ a cardinal number. Recall that $d(X)$ stands for the density of $X$.

\begin{proposition}\label{prop_caliber_density} If $d(X)<\cf(\kappa)$, then $\kappa$ is a caliber for $X$. In particular, $\{\lambda\in\CN : \cf(\lambda)>d(X)\}\subseteq \C(X)$.

\end{proposition}

\begin{proposition}\label{prop_weak_precaliber_celullarity} If $\kappa$ is a weak precaliber for $X$ and $\mathcal{U}$ is a cellular family in $X$, then $|\mathcal{U}|<\kappa$.

\end{proposition}

\begin{proposition}\label{prop_caliber_cofinality} If $\kappa$ a caliber (resp., precaliber; weak precaliber) for $X$, then so is $\cf(\kappa)$.

\end{proposition}

\begin{proposition}\label{prop_weak_precaliber_Hausdorff} If $X$ is Hausdorff, then $\WP(X)\subseteq \UC$.

\end{proposition}

The following result can be found in \cite[Theorem~3.37, p.~14]{riotam2023}.

\begin{theorem}\label{thm_caliber_countable} If $X$ is countably infinite and $T_1$, then $\C(X)=\UC$.

\end{theorem}

\section{Chain conditions in Pixley-Roy hyperspaces}\label{sec_pixroy}

For a topological space $X$, $\F(X)$ is the family $[X]^{<\omega}\setminus\{\emptyset\}$. For every $F\in \F(X)$ and $S\subseteq X$ define $$[F,S] := \{G\in \F(X) : F\subseteq G \subseteq S\}.$$ The family $\{[F,U] : F\in \F(X) \ \wedge \ U\in \tau_X\}$ is a basis for a topology on $\F[X]$ known as the \textit{Pixley-Roy topology}. To follow the traditional notation of the literature, the symbol $\F[X]$ will stand for the set $\F(X)$ equipped with the Pixley-Roy.

These spaces were introduced in 1969 by Pixley and Roy in \cite{pixroy1969} with the aim of presenting a non-separable Moore space with countable cellularity. It is well known that if $X$ is a $T_1$-space, then $\F[X]$ is zero-dimensional (hence, completely regular), Hausdorff and hereditarily metacompact (see \cite{vanDouwen1977}). 

We are interested in determining precisely who the collections $\C\left(\F[X]\right)$, $\P\left(\F[X]\right)$, and $\WP\left(\ F[X]\right)$. Notice that when $X$ is finite then so is $\F[X]$; thus, all infinite cardinals are calibers for $\F[X]$. For this reason, we must enforce the constraint $|X|\in\CN$. Also, we will require $X$ to be a $T_1$-space so that the $\F[X]$ has nice separation properties.

With these conventions in mind, let us start by recalling a basic result regarding the density and cellularity of $\F[X]$ (see \cite[Theorem~2, p.~337]{sakai1983}).

\begin{proposition}\label{prop_pixroy_cardinal_functions} If $X$ is a topological space, then $d\left(\F[X]\right)=|X|$ and $hd(X)\cdot hL(X) \leq c\left(\F[X]\right)\leq nw(X)$.

\end{proposition}

\begin{remark}\label{obs_F[S]} It is appropriate to mention that in Proposition~\ref{prop_pixroy_cardinal_functions} equality is not necessarily reached in the relation $hd(X)\cdot hL(X) \leq c\left(\F[X]\right)$. For example, if we denote the Sorgenfrey line by $\mathbb{S}$ and for each $x\in \mathbb{S}$ we define $U_x: = [x,\infty)$, then it is easy to check that the collection $\left\{[\{x\},U_x] : x\in \mathbb{S}\right\}$ is a cellular family in $\F[\mathbb{S}]$ of size $\mathfrak {c}$, while $hd(\mathbb{S})\cdot hL(\mathbb{S}) = \omega$.

\end{remark}

Now, the $\F[X]$ hyperspaces have a particular feature in that their collection of calibers can be easily determined, but their families of precalibers and weak precalibers are hard to find. Let us first focus on the collection $\C\left(\F[X]\right)$.

\begin{theorem}\label{thm_pixroy_caliber} If $X$ is a topological space, then $$\C\left(\F[X]\right)=\left\{\kappa\in\CN : \cf(\kappa)>|X|\right\}.$$

\end{theorem}

\begin{proof} Combining Propositions~\ref{prop_caliber_density} and \ref{prop_pixroy_cardinal_functions} we immediately deduce that $\{\kappa\in\CN : \cf(\kappa)>|X|\}\subseteq \C\left(\F[X]\right)$. On the other hand, if $\kappa\in\CN$ is such that $\cf(\kappa)\leq |X|$, then we fix a subset $\{x_\alpha : \alpha<\cf(\kappa) \}$ of $X$ that is enumerated without repetitions. Thus, there is no such thing as $J\in [\cf(\kappa)]^{\cf(\kappa)}$ with $\bigcap\left\{[\{x_\alpha\} ,X] : \alpha\in J\right\}\neq \emptyset$ since, if $F\in \bigcap\left\{[\{x_\alpha\} ,X] : \alpha\in J\right\}$, then $F$ would be a finite subset of $X$ containing the infinite set $\{x_\alpha : \alpha\in J\}$, which is absurd. This shows that $\cf(\kappa)$ is not a caliber for $\F[X]$ and hence we deduce from Proposition~\ref{prop_caliber_cofinality} that $\kappa$ is also not a caliber for $\F[X]$.
\end{proof}

Thus, a combination of Proposition~\ref{prop_weak_precaliber_Hausdorff} with Theorem~\ref{thm_pixroy_caliber} implies that if $X$ is a topological space, then \begin{align}\label{pixroy_inclusion} \left\{\kappa\in\CN : \cf(\kappa)>|X|\right\}\subseteq \P\left(\F[X]\right)\subseteq \WP\left(\F[X]\right) \subseteq \UC.
\end{align}

In particular, Theorem~\ref{thm_caliber_countable} and the relations in (\ref{pixroy_inclusion}) imply the following result.

\begin{corollary} If $X$ is a countably infinite space, then \[\UC=\C(X)=\C\left(\F[X]\right)=\P\left(\F[X]\right)=\WP\left(\F[X]\right).\] 

\end{corollary}

Furthermore, it is well known that many cardinal functions of $\F[X]$ coincide with $|X|$ (see \cite[Theorem~2, p.~337]{sakai1983}). In the same vein as these results, if we define the {\it \v{S}anin number} of $X$ as the cardinal $$\check{s}(X) := \min\left\{\kappa \in\CN : \kappa^+ \ \text{is a caliber for} \ X\right\},$$ then Theorem~\ref{thm_pixroy_caliber} implies that $\check{s}\left(\F[ X]\right) = |X|$, in other words, the cardinal number $\check{s}\left(\F[X]\right)$ also enters the list of cardinal functions of $\F[X ]$ that match the cardinality of $X$.

To calculate precisely who the collections $\P\left(\F[X]\right)$ and $\WP\left(\F[X]\right)$ are, we present the following definition.

\begin{definition}\label{def_C(kappa)} If $X$ is a topological space and $\kappa$ is a cardinal number, we say that $X$ satisfies the condition $C(\kappa)$ if and only if for any $\{x_\alpha : \alpha<\kappa\}\subseteq X$ and $\{U_\alpha : \alpha<\kappa\}\subseteq \tau_X$ with $x_\alpha \in U_\alpha$ provided that $\alpha<\kappa$, there exists $J\in [\kappa]^{\kappa}$ such that $\{x_\alpha : \alpha\in J\} \subseteq \bigcap\{U_\alpha : \alpha\in J\} $.

\end{definition}

\begin{remark}\label{rmk_C(kappa)} Every topological space that satisfies the condition $C(\kappa)$ necessarily has caliber $\kappa$. Furthermore, it is clear that if $Y$ is a subspace of $X$ and $X$ satisfies the property $C(\kappa)$, then $Y$ does too; that is, the property $C(\kappa)$ is hereditary.
\end{remark}

We will show later in Theorem~\ref{thm_C(kappa)_equivalence} what is the relationship between $C(\kappa)$ and the chain conditions of Pixley-Roy hyperspaces. What follows is to establish general properties about the $C(\kappa)$ condition that we will use in the rest of this section. We will now show some results regarding the preservation of the condition $C(\kappa)$ in certain topological constructions. 

Before we begin, it is worth mentioning two things: first, the behavior of the $C(\kappa)$ property is similar to that of usual chain conditions; and second, from Lemma~\ref{lemma_C(kappa)_continuous_image} to Proposition~\ref{prop_C(kappa)_sum} it is not necessary that the spaces considered be $T_1$.

The following lemma can be proven by taking preimages and choosing points adequately.

\begin{lemma}\label{lemma_C(kappa)_continuous_image} Let $\kappa$ be a cardinal number and $f: X \to Y$ a continuous surjective function between topological spaces. If $X$ satisfies $C(\kappa)$, then so does $Y$.

\end{lemma}

\begin{proposition}\label{prop_C(kappa)_cofinality} Let $\kappa$ be a cardinal number and $X$ a topological space. If $X$ satisfies $C(\kappa)$, then $X$ also satisfies $C(\cf(\kappa))$.

\end{proposition}

\begin{proof} Clearly the statement is true when $\kappa$ is regular. Suppose then that $\kappa$ is singular and fix a strictly increasing sequence $\{\kappa_\alpha : \alpha<\cf(\kappa)\}$ formed by cardinal numbers in such a way that $\kappa_0 = 0$, $\sup\{\kappa_\alpha : \alpha<\cf(\kappa)\} = \kappa$ and if $\alpha<\cf(\kappa)$ is a limit ordinal, then $\kappa_\alpha = \sup\{\kappa_\beta : \beta<\alpha\}$. Also, let $\{U_\alpha : \alpha<\cf(\kappa)\} \subseteq \tau_X$ and $\{x_\alpha : \alpha<\cf(\kappa)\}$ satisfy $x_\alpha \in U_\alpha$ provided that $\alpha<\cf(\kappa)$. Let us define collections $\{V_\beta : \beta<\kappa\} \subseteq \tau_X$ and $\{y_\beta : \beta<\kappa\} \subseteq X$ as follows: $V_\beta = U_\gamma$ and $y_\beta = x_\gamma$ if and only if $\kappa_\gamma$ is the only cardinal in $\{\kappa_\alpha : \alpha<\cf(\kappa)\}$ with $\kappa_\gamma \leq \beta \leq \kappa_{\gamma+1}$. Then, since $X$ satisfies $C(\kappa)$ there exists $I\in [\kappa]^{\kappa}$ with $\{y_\beta : \beta \in I\} \subseteq \bigcap\{V_\beta : \beta \in I\}$. Thus, if $J:= \{\alpha < \cf(\kappa) : \exists \beta \in I( V_\beta = U_\alpha)\}$, then $J = \{\alpha < \cf (\kappa) : [\kappa_\alpha, \kappa_{\alpha+1})\cap I\neq \emptyset\}$ and hence $|J|=\cf(\kappa)$. Finally, $J \in[\cf(\kappa)]^{\cf(\kappa)}$ and $\{x_\alpha : \alpha \in J\} \subseteq \bigcap\{U_\alpha : \alpha \in J\}$.
\end{proof}

\begin{lemma}\label{lemma_C(kappa)_union} Let $\kappa$ and $\lambda$ be a pair of cardinal numbers ($\lambda$ not necessarily infinite), $X$ a topological space and $\{X_\alpha : \alpha<\lambda\}$ a family of subspaces of $X$. If $\cf(\kappa)>\lambda$ and each $X_\alpha$ satisfies $C(\kappa)$, then so does $\bigcup\{X_\alpha : \alpha<\lambda\}$.

\end{lemma}

\begin{proof} Let $Y:=\bigcup\{X_\gamma : \gamma<\lambda\}$, $\{x_\gamma : \gamma<\kappa\}\subseteq Y$ and $\{U_\gamma : \gamma <\kappa\}\subseteq \tau_Y^+$ with $x_\gamma \in U_\gamma$ for each $\gamma<\kappa$. For every $\gamma<\kappa$ let $\alpha(\gamma)<\lambda$ be such that $x_\gamma \in X_{\alpha(\gamma)}$. The inequality $\cf(\kappa)>\lambda$ ensures that the function $\kappa \to \lambda$ defined as $\gamma\mapsto \alpha(\gamma)$ admits a fiber of size $\kappa$; that is, there exist $J\in [\kappa]^{\kappa}$ and $\alpha<\lambda$ such that $\alpha(\gamma)=\alpha$, provided that $\gamma\in J$. Lastly, since $\{U_\gamma\cap X_\alpha : \gamma\in J\}$ is a subset of $\tau_{X_\alpha}^{+}$ and $X_\alpha$ satisfies $C (\kappa)$, we get $I\in [J]^{\kappa}$ with $\{x_\alpha : \alpha\in I\}\subseteq \bigcap\{U_\gamma\cap X_\alpha : \gamma\in I\}$; consequently, $Y$ satisfies $C(\kappa)$.
\end{proof}

\begin{proposition}\label{prop_C(kappa)_sum} Let $\kappa$ and $\lambda$ be a pair of cardinal numbers ($\lambda$ not necessarily infinite). If $\{X_\alpha : \alpha<\lambda\}$ is a family of topological spaces, then $\bigoplus\{X_\alpha : \alpha<\lambda \}$ satisfies $C(\kappa)$ if and only if $\cf(\kappa)>\lambda$ and each $X_\alpha$ satisfies $C(\kappa)$.

\end{proposition}

\begin{proof} The reciprocal implication follows from Lemma~\ref{lemma_C(kappa)_union}. To check the direct implication, suppose that $X:=\bigoplus\{X_\alpha : \alpha<\lambda \}$ satisfies $C(\kappa)$. Since each $X_\alpha$ is a subspace of $X$, Remark~\ref{rmk_C(kappa)} implies that $X_\alpha$ satisfies $C(\kappa)$. On the other hand, if we assume that $\cf(\kappa)\leq \lambda$, then $\cf(\kappa)$ is not a weak precaliber for $X$ since $\{X_\alpha : \alpha<\cf(\kappa)\}$ is a cellular family and thus, Proposition~\ref{prop_caliber_cofinality} ensures that $\kappa$ is also not a weak precaliber for $X$, a contradiction to our hypothesis (see Remark~\ref{rmk_C(kappa)}). In sum, $\cf(\kappa)>\lambda$.
\end{proof}

In Section~\ref{sec_C(kappa)_products} we will determine the interaction between the condition $C(\kappa)$ and the topological products. For now, it is a good idea to start making connections between the $C(\kappa)$ property and Pixley-Roy hyperspaces with the results we have shown so far.

In \cite[Theorem~3.18, p.~3087]{sakai2012} Sakai proved that if $X$ is a topological space, then $\F[X]$ has precaliber $\omega_1$ if and only if $X$ satisfies $C(\omega_1)$. As we show below, Sakai's argument can be adapted for any cardinal number with uncountable cofinality to obtain a fundamental equivalence for this section of the text.

\begin{theorem}\label{thm_C(kappa)_equivalence} If $X$ is a topological space and $\kappa$ is a cardinal number with $\cf(\kappa)>\omega$, then the following statements are equivalent.

\begin{enumerate}
\item $\F[X]$ has precaliber $\kappa$.
\item $\F[X]$ has weak precaliber $\kappa$.
\item $X$ satisfies $C(\kappa)$. 
\end{enumerate}

\end{theorem}

\begin{proof} Clearly (1) implies (2). To verify that (3) follows from (2) suppose $\{x_\alpha : \alpha<\kappa\}\subseteq X$ and $\{U_\alpha : \alpha<\kappa\}\subseteq \tau_X^+$ are such that $x_\alpha \in U_\alpha$, provided that $\alpha<\kappa$. Then, since $\left\{[\{x_\alpha\}, U_\alpha] : \alpha<\kappa \right\}$ is a subset of $\tau_{\F[X]}^{+}$ and $\F[X]$ has weak precaliber $\kappa$, we deduce the existence of $J\in [\kappa]^{\kappa}$ such that $\left\{[\{x_\alpha\}, U_\alpha] : \alpha\in J \right\}$ is linked. Finally, if $\alpha,\beta \in J$ the relation $[\{x_\alpha\}, U_\alpha] \cap [\{x_\beta\}, U_\beta] \neq \emptyset$ implies that $x_\alpha \in U_\beta$; consequently, $\{x_\alpha : \alpha\in J\} \subseteq \bigcap\{U_\alpha : \alpha\in J\}$.

To show that (3) implies (1) let us set up a collection $\{[F_\alpha, U_\alpha] : \alpha<\kappa\}$ of basic, non-empty open subsets of $\F[X]$. Since $\cf(\kappa)>\omega$ we can assume, without loss of generality, that there exists $n\in\mathbb{N}$ such that $|F_\alpha|=n$ for each $\alpha<\kappa$. Let $\{x(\alpha,k) : k<n\}$ be an enumeration without repetitions of $F_\alpha$ for all $\alpha<\kappa$.

Recursively construct a collection $\{J_k : k<n\}$ such that the following conditions are satisfied for any $k<n$:

\begin{enumerate}
\item $J_0 \in [\kappa]^{\kappa}$;
\item if $k+1<n$, $J_{k+1} \in [J_k]^{\kappa}$; and
\item $\{x(\alpha,k) : \alpha\in J_k\}\subseteq \bigcap\{U_\alpha : \alpha\in J_k\}$.
\end{enumerate}

Thus, $J_{n-1} \in [\kappa]^{\kappa}$ and $\bigcup\left\{F_\alpha : \alpha\in J_{n-1}\right\} \subseteq \bigcap\{U_\alpha : \alpha\in J_{n-1}\}$. Finally, if $I\in [J_{n-1}]^{<\omega}$, then $F:=\bigcup\{F_\alpha : \alpha\in I\}$ satisfies $F\in \bigcap \left\{[F_\alpha, U_\alpha] : \alpha\in I \right\}$; in other words, $\left\{[F_\alpha, U_\alpha] : \alpha\in J_{n-1} \right\}$ is a centered family.
\end{proof}

\begin{remark}\label{rmk_C(kappa)_far} A consequence of Remark~\ref{rmk_C(kappa)} is that if a space $X$ satisfies the condition $C(\kappa)$, then $X$ has caliber $\kappa$ hereditarily. For this reason, Theorem~\ref{thm_C(kappa)_equivalence} implies that the last inclusion of (\ref{pixroy_inclusion}) can be strengthened to $$\WP\left(\F[X]\right) \subseteq \{\kappa \in \CN : X \ \text{has caliber} \ \kappa \ \text{hereditarily} \}\cap \UC.$$

However, in general, this relation is not an equality. For example, the Sorgenfrey line $\mathbb{S}$ is hereditarily separable (thus, it has hereditary caliber $\omega_1$), but $\F[\mathbb{S}]$ admits a cellular family of cardinality $\mathfrak{c}$ (see Remark~\ref{obs_F[S]}); in particular, $\F[\mathbb{S}]$ does not have weak precaliber $\omega_1$.

To give us an idea of the strength of the condition $C(\omega_1)$, Sakai mentions in \cite[Corollary~3.19, p.~3087]{sakai2012} that if a space $X$ satisfies $C(\omega_1 )$, then necessarily $X^{\omega}$ is hereditarily separable and hereditarily Lindel\"of.

\end{remark}

The last containment of (\ref{pixroy_inclusion}) and Theorem~\ref{thm_C(kappa)_equivalence} imply the following result.

\begin{corollary}\label{cor_thm_C(kappa)_equivalence} If $X$ is a topological space, then \[\P\left(\F[X]\right)=\WP\left(\F[X]\right)=\left\{\kappa\in\CN :  X \ \text{satisface la condición} \ C(\kappa)\right\}\cap \UC.\]

\end{corollary}

In light of Corollary~\ref{cor_thm_C(kappa)_equivalence}, our next objective is to try to find internal conditions to know when a space satisfies or not the condition $C(\kappa)$ for a cardinal $\kappa$ with uncountable cofinality.

First, we show below that when we consider the net weight of the space in a convenient way, we obtain a positive answer.

\begin{proposition}\label{prop_C(kappa)_nw} If $X$ is a topological space and $\kappa$ is a cardinal number such that $nw(X)<\cf(\kappa)$, then $X$ satisfies $C(\kappa)$. In particular, $\{\kappa\in\CN : \cf(\kappa)>nw(X)\} \subseteq \P\left(\F[X]\right)$.

\end{proposition}

\begin{proof} Let $\mathcal{N}$ be a net for $X$ of minimum cardinality, and $\{x_\alpha : \alpha<\kappa\}\subseteq X$ and $\{U_\alpha : \alpha<\kappa\}\subseteq \tau_X^+$ such that $x_\alpha \in U_\alpha$, provided that $\alpha<\kappa$. For each $\alpha<\kappa$, let $N_\alpha \in \mathcal{N}$ satisfy $x_\alpha \in N_\alpha \subseteq U_\alpha$. It turns out that the function $\kappa \to \mathcal{N}$ given by $\alpha \mapsto N_\alpha$ has a fiber of cardinality $\kappa$; that is, there exists $J\in [\kappa]^{\kappa}$ such that $N_\alpha=N_\beta$ whenever $\alpha,\beta \in J$. Thus, $\{x_\alpha : \alpha\in J\} \subseteq \bigcap\{U_\alpha : \alpha\in J\} $.
\end{proof}

\begin{corollary} If $X$ is a cosmic space, then \[\P\left(\F[X]\right)=\UC=\WP\left(\F[X]\right).\] 

\end{corollary}

An immediate consequence of Theorem~\ref{thm_pixroy_caliber} and Proposition~\ref{prop_C(kappa)_nw} is that if $X$ satisfies $nw(X)<|X|$, then the collections $\C\left (\F[X]\right)$ and $\P\left(\F[X]\right)$ do not match due to the relations \begin{align*} \C\left(\F[X]\right) &= \left\{\kappa\in\CN : \cf(\kappa)>|X|\right\} \\
&\subsetneq \left\{\kappa\in\CN : \cf(\kappa)>nw(X)\right\} \subseteq \P\left(\F[X]\right).
\end{align*} Notice that $nw(X)^+ \in \P\left(\F[X]\right) \setminus \C\left(\F[X]\right)$. In sum, the following corollary is verified.

\begin{corollary} If $X$ is a topological space with $nw(X)<|X|$, then $\C\left(\F[X]\right) \subsetneq \P\left(\F[X]\right)$.

\end{corollary}

Now, our next results are intended to set conditions on $X$ to ensure that the condition $C(\kappa)$ is not satisfied on $X$.

\begin{proposition}\label{prop_C(kappa)_d} If $X$ is a topological space and $\kappa$ is a cardinal number such that $\kappa\leq d(X)$, then $X$ does not satisfy $C(\kappa)$.

\end{proposition}

\begin{proof} Our goal is to recursively construct a pair of collections $\{x_\alpha : \alpha<\kappa\}\subseteq X$ and $\{U_\alpha : \alpha<\kappa\}\subseteq \tau_X^+$ such that, for each $\alpha<\kappa$, $x_\alpha \in U_\alpha$ and $U_\alpha \cap \{x_\beta : \beta<\alpha\}=\emptyset$. Suppose that for $\alpha<\kappa$ we have found families $\{x_\beta : \beta<\alpha\}\subseteq X$ and $\{U_\beta : \beta<\alpha\}\subseteq \tau_X ^+$ with the desired conditions. Then, since $\{x_\beta : \beta<\alpha\}$ is not a dense subset of $X$, there exists $U_\alpha \in \tau_X^+$ such that $U_\alpha \cap \{ x_\beta : \beta<\alpha\}=\emptyset$. Thus, if $x_\alpha \in U_\alpha$, then $\{x_\beta : \beta<\alpha+1\}$ and $\{U_\beta : \beta<\alpha+1\}$ have the required properties.

From the way the previous sets were constructed, it is clear that for any $\alpha<\kappa$ the relation $x_\alpha \not \in \bigcup\{U_\beta : \beta>\alpha\}$ is satisfied. Finally, if $J\in [\kappa]^{\kappa}$ there exists $\beta \in J \cap (\alpha,\kappa)$ and hence, $x_\alpha \not \in U_ \beta$; consequently, $\{x_\alpha : \alpha\in J\} \not\subseteq \bigcap\{U_\alpha : \alpha\in J\}$.
\end{proof}

\begin{corollary}\label{cor_C(kappa)_hd} If $X$ is a topological space and $\kappa$ is a cardinal number such that $\cf(\kappa)< hd(X)$, then $X$ does not satisfy $C(\kappa)$.

\end{corollary}

\begin{proof} The relation $\cf(\kappa)< hd(X)$ implies the existence of $Y\subseteq X$ with $\cf(\kappa)<d(Y)$. Then, by Proposition~\ref{prop_C(kappa)_d} it is verified that $Y$ does not satisfy $C(\cf(\kappa))$, and Proposition~\ref{prop_C(kappa)_cofinality} ensures that $Y$ does not satisfy $C(\kappa)$. Lastly, since $C(\kappa)$ is hereditary (see Remark~\ref{rmk_C(kappa)}), $X$ does not satisfy $C(\kappa)$.
\end{proof}

Thanks to Proposition~\ref{prop_C(kappa)_nw} and Corollary~\ref{cor_C(kappa)_hd}, the only cardinal numbers $\kappa$ that remain to be analyzed are those that have uncountable cofinality and satisfy the relations \begin{align}\label{pixroy_hd_nw} hd(X)\leq \cf(\kappa)\leq nw(X).
\end{align}

To refine the first inequality of (\ref{pixroy_hd_nw}) it is necessary to further expand our conceptual field.

\section{Weakly separated spaces}\label{sec_weak_sep}

Unless explicitly stated otherwise, the topological spaces in this section will not be constrained to satisfy the separation axiom $T_1$.

Recall that if $X$ is a set, $x\in X$ and $<$ is an ordering relation on $X$, then the {\it initial segment} determined by $x$ is the set $(\gets ,x) := \{y\in X : y<x\}$.

A topological space $X$ is {\it right-separated} (resp., {\it left-separated}) if it admits a well order such that its initial segments are elements of $\tau_X$ (resp., of $\tau_X^{*}$). The {\it height} and {\it width} of $X$ are, respectively, the cardinal numbers\begin{align*}  h(X) &:= \sup \{|Y| : Y \ \text{is a right-separated subset of} \ X\}+\omega \ \text{and} \\
z(X) &:= \sup \{|Y| : Y \ \text{is a left-separated subset of} \ X\}+\omega.
\end{align*}

These types of spaces have been extensively studied in the literature, even in connection with the chain conditions of topological spaces (see \cite{juhasz1980} and \cite{juhshe2003}). A well known result is that $hL = h$ and $hd = z$ (see \cite[2.9, p.~16]{juhasz1980}). 

On the other hand, a space $X$ is {\it weakly separated} if there exists a family $\{U_x : x\in X\} \subseteq \tau_X$ that satisfies the following conditions: $x\in U_x$ for each $x\in X$, and if $x,y\in X$ are distinct, then $x\not\in U_y$ or $y\not \in U_x$. These spaces were introduced by Tkachenko in \cite{tkachenko1978}. In these circumstances we will say that the family $\{U_x : x\in X\}$ is a {\it weak separation} for $X$. Finally, the {\it weak separation number} of $X$ is the cardinal number $$R(X) := \sup\{|Y| : Y \ \text{is a weakly separated subspace of} \ X\}+\omega.$$

It is easy to construct weakly separated spaces, e.g. any countable $T_1$-space has this property. Indeed, if $X$ is finite then $\{\{x\} : x\in X\}$ is a weak separation for $X$ (any discrete space is weakly separated). On the other hand, if $\{x_n : n<\omega\}$ is an enumeration without repetitions of $X$ and for each $n<\omega$ we define $U_n:= X \setminus \{x_k : k<n \}$, then $\{U_n : n<\omega\}$ is a weak separation for $X$.

One more example that is fundamental to us is the Sorgenfrey line. Just notice that the family $\{[x,\infty) : x\in \mathbb{S}\}$ is a weak separation for $\mathbb{S}$.

Also, it is not difficult to check that all right-separated or left-separated subspaces of a space $X$ are weakly separated; consequently, $hd(X)\cdot hL(X) \leq R(X)$. Furthermore, as for any weakly separated space it is true that $|\cdot|=nw$, then $R(X) \leq nw(X)$. We collect the observations of this paragraph in the following result.

\begin{proposition}\label{prop_weak_sep_car_fun} If $X$ is a topological space $hd(X)\cdot hL(X) \leq R(X) \leq nw(X)$.

\end{proposition}

\begin{corollary}\label{cor_weak_sep_uncountable} If $nw(X)<|X|$, then $X$ is not weakly separated. In particular, uncountable cosmic spaces are not weakly separated.

\end{corollary}

We now expose some basic properties of weakly separated spaces. For example, routine arguments can be used to prove the following proposition.

\begin{proposition}\label{prop_weak_sep_basic} Let $X$, $Y$ and $Z$ be a triplet of topological spaces.

\begin{enumerate}
\item If $Y$ is a subspace of $X$ and $Z$ is a subspace of $Y$, then $Z$ is weakly separated as a subspace of $X$ if and only if $Z$ is weakly separated as a subspace of $Y$.
\item When $X$ is weakly separated, $X$ is $T_0$.
\item If $X$ is weakly separated and $Y$ is a subspace of $X$, then $Y$ is weakly separated.
\item The following statements are true for a function $f:X\to Y$.

\begin{enumerate}
\item If $f$ is a condensation and $Y$ is weakly separated, then $X$ is weakly separated.

\item If $f$ is continuous and surjective, then $R(Y)\leq R(X)$.
\end{enumerate}
\end{enumerate}

\end{proposition}

\begin{proposition}\label{prop_weak_sep_union_open} If $X$ is a topological space and $\mathcal{U}$ is a family of weakly separated open subspaces of $X$, then $\bigcup\mathcal{U}$ is a weakly separated subspace of $X$.

\end{proposition}

\begin{proof} Let $U:= \bigcup\mathcal{U}$, $\{U_\alpha : \alpha<\lambda\}$ be an enumeration without repetitions of $\mathcal{U}$ and, for each $\alpha<\lambda$, fix a weak separation $\{V(x,\alpha) : x\in U_\alpha\}\subseteq \tau_{U_\alpha}$ for $U_\alpha$. Furthermore, for each $x\in U$ consider the ordinal $\alpha_x := \min\{\alpha<\lambda : x\in U_\alpha\}$ and define $V_x := V(x,\alpha_x) $. Now, to check that $\{V_x : x\in U\}$ is a weak separation for $U$, first note that $\{V_x : x\in U\}$ is a subset of $\tau_U$ with $x\in V_x$ for all $x\in U$. Finally, if $x,y\in U$ are different we have two cases: if $\alpha_x < \alpha_y$ or $\alpha_y < \alpha_x$, then $y \not \in V_x$ or $x \not \in V_y$; otherwise, $\alpha_x = \alpha_y$ and so, since $U_{\alpha_x}$ is weakly separated, $x \not \in V_y$ or $y \not \in V_x$.
\end{proof}

\begin{proposition} Let $X$ be a topological space and $\mathcal{F}$ a family of weakly separated closed subspaces of $X$. If $\mathcal{F}$ is locally finite with respect to $\bigcup\mathcal{F}$, then $\bigcup\mathcal{F}$ is a weakly separated subspace of $X$.

\end{proposition}

\begin{proof} Let $F:= \bigcup\mathcal{F}$ and $\{F_\alpha : \alpha<\lambda\}$ be an enumeration without repetitions of $\mathcal{F}$. For each $x\in F$ let $U_x \in \tau_{F}$ be such that $x\in U_x$ and $I_x := \{\alpha<\lambda : F_\alpha \cap U_x \neq \emptyset\}$ is a finite set, and define $J_x := \{\alpha \in I_x : x\in F_\alpha\}$ and $K_x:= I_x \setminus J_x$. Furthermore, for all $\alpha<\lambda$ let $\{V(x,\alpha) : x\in F_\alpha\}$ be a subset of $\tau_F$ that weakly separates $F_\alpha$. Finally, for each $x\in F$ define $$W_x := U_x \cap \left(\bigcap_{\alpha \in J_x} V(x,\alpha)\right) \cap \left(\bigcap_{\alpha \in K_x} (X\setminus F_\alpha)\right).$$

We will show $\{W_x : x\in F\}$ is a weak separation for $F$. Note that $\{W_x : x\in F\}$ is a subset of $\tau_F$ with $x\in W_x$ for all $x\in F$. Now, let $x,y\in F$ be distinct and $\alpha,\beta<\lambda$ be such that $x\in F_\alpha$ and $y\in F_\beta$. Suppose further that $x\in W_y$ and observe that, since $x\in F_\alpha\cap U_y$, $\alpha$ belongs to $I_y$. Then, since for each $\gamma \in K_y$ we have that $x\in X\setminus F_\gamma$, we infer that $\alpha \not \in K_y$ and, therefore, $\alpha\in J_y$ . Thus, $x$ and $y$ are elements of $F_\alpha$ such that $x\in V(y,\alpha)$; consequently, since $\{V(z,\alpha) : z\in F_\alpha\}$ is a weak separation for $F_\alpha$, we obtain the relation $y\not \in V(x,\alpha)$, which implies that $y\not \in W_x$.
\end{proof}

The direct implication of the following result follows from Proposition~\ref{prop_weak_sep_basic}(3), while the converse implication is a consequence of Proposition~\ref{prop_weak_sep_union_open}.

\begin{proposition} If $\lambda$ is a cardinal number (not necessarily infinite) and $\{X_\alpha : \alpha<\lambda\}$ is a family of topological spaces, then $\bigoplus_{\alpha<\lambda} X_\alpha$ is weakly separated if and only if each $X_\alpha$ is weakly separated.

\end{proposition}

\begin{proposition}\label{prop_weak_sep_sum_R} If $\lambda$ is a cardinal number (not necessarily infinite) and $\{X_\alpha : \alpha<\lambda\}$ is a family of topological spaces, then $R\left(\bigoplus_{\alpha<\lambda} X_\alpha\right) = \lambda \cdot \sup\left\{R(X_\alpha) : \alpha<\lambda\right\}$.

\end{proposition}

\begin{proof} Let $X:= \bigoplus_{\alpha<\lambda} X_\alpha$. First, if for each $\alpha<\lambda$ we take $x_\alpha \in X_\alpha$, then $\{x_\alpha : \alpha<\lambda\}$ is a discrete subspace (in particular, weakly separated) of $X$ and therefore $\lambda \leq R(X)$. On the other hand, if $\alpha<\lambda$ and $Y_\alpha$ is a weakly separated subspace of $X_\alpha$, then $Y_\alpha$ is a weakly separated subspace of $X$ (see proposition ~ \ref{prop_weak_sep_basic}); consequently, $\sup\left\{R(X_\alpha) : \alpha<\lambda\right\} \leq R(X)$.

To check the remaining inequality, let $Y$ be a weakly separated subspace of $X$. Since for every $\alpha<\lambda$ it is satisfied that $Y \cap X_\alpha$ is a weakly separated subspace of $X_\alpha$ (see proposition~\ref{prop_weak_sep_basic}), we deduce that $|Y \cap X_\alpha|\leq R(X_\alpha)$; thus $$|Y| = \left|\bigcup_{\alpha<\lambda} Y \cap X_\alpha\right| \leq \lambda\cdot \sup\{|Y \cap X_\alpha| : \alpha<\lambda\} \leq \lambda \cdot \sup\left\{R(X_\alpha) : \alpha<\lambda\right\}.$$ In conclusion, $R(X) \leq \lambda \cdot \sup\left\{R(X_\alpha) : \alpha<\lambda\right\}$.
\end{proof}

\begin{proposition}\label{prop_weak_sep_product} If $\lambda$ is a cardinal number (not necessarily infinite) and $\{X_\alpha : \alpha<\lambda\}$ is a family of topological spaces with al least two points, then the following statements are true.

\begin{enumerate}
\item If $\prod_{\alpha<\lambda} X_\alpha$ is weakly separated, then each $X_\alpha$ is weakly separated.
\item Whenever each $X_\alpha$ is weakly separated, the box product $\square_{\alpha<\lambda} X_\alpha$ is weakly separated. In particular, the topological product of a finite family of weakly separated spaces is weakly separated.
\item $\prod_{\alpha<\lambda} X_\alpha$ is not weakly separated if $\lambda\geq \omega$.

\end{enumerate}

\end{proposition}

\begin{proof} For part (1) it is enough to remember that each $X_\alpha$ is homeomorphic to a subspace of the topological product $\prod_{\alpha<\lambda} X_\alpha$.

Second, for each $\alpha<\lambda$ let $\{U(\alpha,x) : x\in X_\alpha\}$ be a weak separation for $X_\alpha$. Observe that if $X:= \square_{\alpha<\lambda} X_\alpha$, $x,y\in X$ are distinct and $\alpha<\lambda$ satisfies $x(\alpha) \neq y(\alpha)$, then $x(\alpha) \not \in U(\alpha,y(\alpha))$ or $y(\alpha) \not \in U(\alpha,x(\alpha ))$ and thus $x\not\in \prod_{\beta<\lambda} U(\beta,y(\beta))$ or $y\not\in \prod_{\beta<\lambda } U(\beta,x(\beta))$. Consequently, $\left\{\prod_{\alpha<\lambda} U(\alpha,x(\alpha)) : x\in X\right\}$ is a weak separation for $X$.

For part (3) notice that if $J$ is a countably infinite subset of $\lambda$, for each $\alpha\in J$ we take $ x_\alpha, y_\alpha \in J$ distinct, and we define $Y_\alpha := \{x_\alpha,y_\alpha\}$, then the product $\prod_{\alpha\in J} Y_\alpha $ is cosmic, uncountable and embeds into $\prod_{\alpha<\lambda} X_\alpha$. Therefore, since Corollary~\ref{cor_weak_sep_uncountable} guarantees that $\prod_{\alpha\in J} Y_\alpha$ is not weakly separated, we conclude that $\prod_{\alpha<\lambda} X_\alpha$ is also not weakly separated.
\end{proof}

By virtue of Proposition~\ref{prop_weak_sep_sum_R}, it would be desirable to obtain a similar formula to calculate the value of $R$ for a product of topological spaces. We have not been able to obtain an equality in the previous sense, but we do have a couple of bounds in the case of $T_2$-spaces.

\begin{proposition}\label{prop_weak_sep_product_R} If $\lambda$ is a cardinal number (not necessarily infinite) and $\{X_\alpha : \alpha<\lambda\}$ is a family of Hausdorff spaces, then $$\lambda \cdot \sup\{R(X_\alpha) : \alpha<\lambda\} \leq R\left(\prod_{\alpha<\lambda} X_\alpha\right) \leq \lambda\cdot \sup\left\{2^{R(X_\alpha) } : \alpha<\lambda\right\}.$$

\end{proposition}

\begin{proof} On the one hand, since $\bigoplus_{\alpha<\lambda} X_\alpha$ embeds into $\prod_{\alpha<\lambda} X_\alpha$, Proposition~\ref{prop_weak_sep_sum_R} implies that $$\lambda \cdot \sup\{R(X_\alpha) : \alpha<\lambda\} \leq R\left(\prod_{\alpha<\lambda} X_\alpha\right).$$ On the other hand, Proposition~\ref{prop_weak_sep_car_fun} guarantees that \begin{align*} R\left(\prod_{\alpha<\lambda} X_\alpha\right) &\leq nw\left(\prod_{\alpha<\lambda} X_\alpha\right) \\
&= \lambda \cdot \sup\{nw(X_\alpha) : \alpha<\lambda\} \leq \lambda\cdot \sup\left\{2^{R(X_\alpha)} : \alpha<\lambda\right\}.
\end{align*}
\end{proof}

As our last basic property regarding weakly separated spaces, the following result is mentioned in \cite{tkachenko1978}.

\begin{proposition}\label{prop_weak_sep_sigma} The $\sigma$-product of a family of weakly separated $T_1$-spaces is weakly separated.

\end{proposition}

Proposition~\ref{prop_weak_sep_sigma} cannot be generalized to $\Sigma$-products. For example, as the $\Sigma$-product $\Sigma(\omega_1) := \{x \in D(2)^{\omega_1} : |\{\alpha<\omega_1 : x(\alpha) \neq 0\}| \leq \omega\}$ contains a homeomorphic copy of $D(2)^{\omega}$ which is cosmic and uncountable, it turns out that $\Sigma(\omega_1)$ is not weakly separated (see Corollary~\ref{cor_weak_sep_uncountable}). 

Now, it is convenient to mention that the inequalities exposed in Proposition~\ref{prop_weak_sep_car_fun} can be strict. For example, the Sorgenfrey line satisfies $hd(\mathbb{S})\cdot hL(\mathbb{S}) = \omega < \mathfrak{c} = R(\mathbb{S})$. The second inequality can also be strict as we will see in Example~\ref{ex_R}.

In sum, Propositions~\ref{prop_pixroy_cardinal_functions} and \ref{prop_weak_sep_car_fun} certify that for any space $X$ the cardinal numbers $R(X)$ and $c\left(\F[X]\right)$ are between $hd(X)\cdot hL(X)$ and $nw(X)$. The natural question that arises is: how are the cardinals $R(X)$ and $c\left(\F[X]\right)$ related?

Before answering the previous question, we are going to introduce a couple more cardinal functions. Following the tradition of the literature, for any topological space $X$ we define $$c^*(X) := \sup \{c(X^{n}) : n\in \mathbb{N}\} \quad \text{and} \quad R^{*}(X) := \sup \{R(X^{n}) : n\in \mathbb{N}\}.$$

The following space appears in \cite{mcintyre2006} and \cite{tkachenko1978}.

\begin{example}\label{ex_R} If $\tilde{\mathbb{R}}$ is the space obtained by equipping the set $\mathbb{R}$ with the topology generated by the subbase $$\mathcal{S} := \{U\setminus S : U\in\tau_\mathbb{R} \ \text{and} \ S \ \text{is a non-trivial convergent sequence in} \ \mathbb{R}\},$$ then $\omega = R^* (\tilde{\mathbb{R}}) < nw(\tilde{\mathbb{R}})$.

\end{example}

The connection between $R(X)$, $c\left(\F[X]\right)$, $R^{*}(X)$, and $c^*\left(\F[X]\right )$ is established in Theorems~\ref{thm_R*_c*_1}, \ref{thm_R*_c*_2} and \ref{thm_R*_c*_3} that we present below.

\begin{theorem}\label{thm_R*_c*_1} If $X$ is a topological space, then $R^{*}(X) = c^*\left(\F[X]\right)$.

\end{theorem}

\begin{proof} Let $n\in\mathbb{N}$, $\{x_\alpha : \alpha<\kappa\}$ be a weakly separated subspace of $X^n$ with $x_\alpha = \left(x(\alpha,1),\dots,x(\alpha,n)\right)$, $\{U(\alpha,1)\times\dots\times U(\alpha,n) : \alpha<\kappa\}$ a weak separation for $\{x_\alpha : \alpha<\kappa\}$, and define for every $\alpha<\kappa$, $U_\alpha := [\{x(\alpha,1)\},U(\alpha,1)]\times \dots \times [\{x(\alpha,n)\},U(\alpha,n)]$. We claim that $\{U_\alpha : \alpha<\kappa\}$ is a cellular family in $\F[X]^n$. Indeed, if $\alpha,\beta<\kappa$ and $(F_1,\dots,F_n)\in \F[X]^n$ are such that $(F_1,\dots,F_n) \in U_\alpha\cap U_\beta$, then for every $1\leq k \leq n$ it is satisfied that $\{x(\alpha,k),x(\beta,k)\} \subseteq F_k \subseteq U(\alpha,k) \cap U(\beta,k)$. Thus, $x_\alpha \in U(\beta,1)\times\dots\times U(\beta,n)$ and $x_\beta \in U(\alpha,1)\times\dots\times U(\alpha,n)$, which implies that $\alpha=\beta$. This argument shows that for any $n\in \mathbb{N}$ the relation $R(X^n) \leq c\left(\F[X]^n\right)$ holds; consequently, $R^*(X) \leq c^*\left(\F[X]\right)$.

For the remaining inequality suppose seeking for a contradiction that $c^*\left(\F[X]\right)$ is strictly greater than $\kappa:= R^*(X)$. Then, since $c^*\left(\F[X]\right)\geq \kappa^+$ there is $n\in \mathbb{N}$ with $c\left(\F[X]^n\right)\geq \kappa^+$ and therefore, there exists a cellular family $\{[F(\alpha,1),U(\alpha,1)]\times \dots \times [F(\alpha,n),U(\alpha,n)] : \alpha<\kappa^+\}$ in $\F[X]^n$. Since the function $\kappa^+ \to \mathbb{N}^{n}$ given by $\alpha \mapsto \left(|F(\alpha,1)|,\dots,|F(\alpha,n)|\right)$ admits a fiber of cardinality $\kappa^+$, we will suppose without loss of generality the existence of $m_1,\dots,m_n \in \mathbb{N}$ in such a way that $|F(\alpha,k)|= m_k$ whenever $\alpha<\kappa^+$ and $1\leq k \leq n$. Let $\{x(\alpha,k,j) : 1\leq j \leq m_k\}$ be an enumeration of $F(\alpha,k)$ with no repetitions.

For every $\alpha<\kappa^{+}$ let \begin{align*} x_\alpha &:= \left(\left(x(\alpha,1,1),\dots,x(\alpha,1,m_1)\right),\dots,\left(x(\alpha,n,1)\dots,x(\alpha,n,m_n)\right)\right) \ \text{and} \\
U_\alpha &:= \left(U(\alpha,1)\times\dots\times U(\alpha,1)\right)\times\dots\times \left(U(\alpha,n)\times\dots\times U(\alpha,n)\right). 
\end{align*} We will prove that $\{U_\alpha : \alpha<\kappa^+\}$ is a weak separation for the subspace $\{x_\alpha : \alpha<\kappa^+\}$ of $X^{m_1}\times\dots\times X^{m_n}$. If $\alpha<\beta<\kappa^+$ and $x_\alpha \in U_\beta$, then for each $1\leq k \leq n$ it is the case that $F(\alpha,k) \subseteq U(\beta,k)$. Hence, there is $1 \leq k \leq n$ such that $F(\beta,k) \not \subseteq U(\alpha,k)$, which implies the existence of $1\leq j \leq n$ with $x(\beta,k,j) \not \in U(\alpha,k)$; consequently, $x_\beta \not \in U_\alpha$. 

Finally, observe that if $m:=m_1+\dots+m_n$, then the spaces $X^{m_1}\times\dots\times X^{m_n}$ and $X^{m}$ are homeomorphic and, therefore, we deduce that $X^{m}$ admits a weakly separated subspace of size $\kappa^+$, a contradiction to the relation $R^*(X) < \kappa^+$.
\end{proof}

\begin{theorem}\label{thm_R*_c*_2} If $X$ is a topological space, then $R(X)\leq c\left(\F[X]\right) \leq R^{*}(X) \leq nw(X)$. If in addition $X$ is $T_2$, then $nw(X) \leq 2^{R(X)}$.

\end{theorem}

\begin{proof} For the inequality $R(X)\leq c\left(\F[X]\right)$ we will show that if $Y$ is a weakly separated subspace of $X$, then $\F[X]$ admits a cellular family of cardinality $|Y|$. Let $\{V_y : y\in Y\}$ be a weak separation for $Y$. For each $y\in Y$ let us take $U_y\in \tau_X$ with $V_y = U_y \cap Y$, and consider the collection $\{[\{y\},U_y] : y\in Y\} $. If $y,z\in Y$ are distinct and we assume that $[\{y\},U_y] \cap [\{z\},U_z] \neq \emptyset$, then $\{y,z\} \subseteq U_y \cap U_z$. In this way, $\{y,z\} \subseteq V_y \cap V_z$ and, therefore, $y\in V_z$ and $z\in V_y$; a contradiction to the weak separation hypothesis. Consequently, $\{[\{y\},U_y] : y\in Y\}$ is a cellular family in $\F[X]$ of size $|Y|$.

The relation $c\left(\F[X]\right) \leq R^{*}(X)$ is evident from Theorem~\ref{thm_R*_c*_1}. On the other hand, since for each $n\in \mathbb{N}$ it is satisfied that $R(X^n) \leq nw(X^n) = nw(X)$ (see Proposition~\ref{prop_weak_sep_car_fun}), then $R^{*}(X) \leq nw(X)$.

The last inequality is a consequence of the following: a remarkable result by Hajnal and Juh\'asz states that $|X| \leq 2^{s(X)\cdot \psi(X)}$ whenever $X$ is $T_1$ (see \cite[Theorem~4.7, p.~20]{hodel1984}). Thus, as for Hausdorff spaces $s\leq hd$ and $\psi \leq hL$, it is immediate that $$nw(X) \leq |X| \leq 2^{s(X)\cdot \psi(X)} \leq 2^{hd(X)\cdot hL(X)}\leq 2^{R(X)}$$ if $X$ is $T_2$ (see Proposition~\ref{prop_weak_sep_car_fun}).
\end{proof}

\begin{corollary}\label{cor_thm_R*_c*_2} Every weakly separated space $X$ satisfies the relations $c^*\left(\F[X]\right)=R^*(X)=c\left(\F[X]\right)=R(X)=nw(X)=|X|$.

\end{corollary}

\begin{theorem}\label{thm_R*_c*_3} If $X$ is a topological space, then $R^{*}(X) = \sup \{c\left(\F[X^n]\right) : n\in \mathbb{N}\}$.

\end{theorem}

\begin{proof} Let $\mu := \sup \{c\left(\F[X^n]\right) : n\in \mathbb{N}\}$. First, since for every $n\in \mathbb{N}$ it is satisfied that $R(X^n)\leq c\left(\F[X^n]\right)$ (see Theorem~\ref{thm_R*_c*_2}), then $R^{*}(X)\leq \mu$. Now suppose for an absurdity that $\mu$ is strictly greater than $\kappa :=R^{*}(X)$. Use the inequality $\mu \geq \kappa^+$ to find $n\in \mathbb{N}$ with $c\left(\F[X^n]\right) \geq \kappa^+$, and let $\{[F_\alpha, U_\alpha] : \alpha<\kappa^+\}$ be a cellular family in $\F[X^n]$. Given that the function $\kappa^+ \to \mathbb{N}$ determined by $\alpha\mapsto |F_\alpha|$ has a fiber of cardinality $\kappa^+$, we will assume without loss of generality that there is $m\in \mathbb{N}$ with $|F_\alpha|=m$ for each $\alpha<\kappa^+$. Let $\{x(\alpha,1),\dots,x(\alpha,m)\}$ be an enumeration without repetitions of $F_\alpha$, and for each $1\leq k \leq m$ let $x (\alpha,k,1),\dots,x(\alpha,k,n) \in X$ and $U(\alpha,k,1),\dots,U(\alpha,k,n) \in \tau_X$ satisfy $x(\alpha,k) = \left(x(\alpha,k,1),\dots,x(\alpha,k,n)\right) \in U(\alpha,k ,1)\times \dots\times U(\alpha,k,n) \subseteq U_\alpha$.

For every $\alpha<\kappa^{+}$ let \begin{align*} x_\alpha &:= \left(\left(x(\alpha,1,1),\dots,x(\alpha,1,n)\right),\dots,\left(x(\alpha,m,1)\dots,x(\alpha,m,n)\right)\right) \ \text{and} \\
U_\alpha &:= \left(U(\alpha,1,1)\times\dots\times U(\alpha,1,n)\right)\times\dots\times \left(U(\alpha,m,1)\times\dots\times U(\alpha,m,n)\right). 
\end{align*} We will see that $\{U_\alpha : \alpha<\kappa^+\}$ is a weak separation for the subspace $\{x_\alpha : \alpha<\kappa^+\}$ of $X^{n}\times\dots\times X^{n}$. If $\alpha<\beta<\kappa^+$ and $x_\alpha \in U_\beta$, then $F_\alpha \subseteq U_\beta$ and therefore, $F_\beta \not\subseteq U_\alpha$, i.e., there are $1\leq k \leq m$ and $1\leq j \leq n$ such that $x(\beta,k,j)\not\in U(\alpha,k,j)$; in particular, $x_\beta \not \in U_\alpha$.

Lastly notice that, since the spaces $X^{n}\times\dots\times X^{n}$ and $X^{mn}$ are homeomorphic, $X^{mn}$ contains a weakly separated subspace of cardinality $\kappa^+$, a contradiction to the inequality $R^*(X) < \kappa^+$.
\end{proof}

\begin{question} Are there examples of spaces $X$ and $Y$ such that $R(X)< c\left(\F[X]\right)$ and $c\left(\F[Y]\right) < R^ *(Y)$?

\end{question}

A consequence of Proposition~\ref{prop_weak_sep_product}(2) is that all finite powers of a weakly separated space also possess the same characteristic. Now, a natural question that might arise is when the hyperspace $\F[X]$ is weakly separated. It turns out that this property is always present. Indeed, simple reasoning shows that $\{[F,X] : F\in \F[X]\}$ is a weak separation for $\F[X]$. Thus, if for each $n\in \mathbb{N}$ we recursively define $\F^1[X] := \F[X]$ and $\F^{n+1}[X] := \F[\F^n[X]]$, then the following result is verified.

\begin{proposition}\label{prop_weak_sep_F[X]_weak_sep} If $X$ is a topological space, then for each $n\in \mathbb{N}$ it is satisfied that $\F^n[X]$, $\F[X^n]$ and $\F[X ]^n$ are weakly separated.

\end{proposition}

In particular, Proposition~\ref{prop_weak_sep_F[X]_weak_sep} allows us to detect by means of the weak separation property when a space $X$ {\bf does not} embed topologically into one of the spaces $\F^n[ X]$, $\F[X^n]$ and $\F[X]^n$.

\begin{corollary}\label{cor_weak_sep_F[X]_weak_sep} If $X$ is not weakly separated, then for each $n\in \mathbb{N}$ it is satisfied that $\F^n[X]$, $\F[X^n]$ and $\F[X ]^n$ does not contain topological copies of $X$.

\end{corollary}

In particular, by Corollary~\ref{cor_weak_sep_uncountable} any space $X$ with $nw(X)<|X|$ satisfies the hypothesis of Corollary~\ref{cor_weak_sep_F[X]_weak_sep}. Thus, for example, $\mathbb{R}$ does not embed into $\F^n[\mathbb{R}]$, $\F[\mathbb{R}^n]$ and $\F[ \mathbb{R}]^n$ for any $n\in \mathbb{N}$.

On the other hand, although there are many difficulties in doing the explicit computation of $c\left(\F[X]\right)$ for an arbitrary topological space $X$, Corollary~\ref{cor_thm_R*_c*_2} and Proposition~\ref{prop_weak_sep_F[X]_weak_sep} allow us to do the corresponding calculation for each $\F^n[X]$ with $n\geq 2$.

\begin{corollary} If $X$ is a topological space and $n\in\mathbb{N}$, then $$c^*\left(\F^{n+1}[X]\right)=R^*\left(\F^{n}[X]\right)=c\left(\F^{n+1}[X]\right)=R\left(\F^{n}[X]\right)=|X|.$$

\end{corollary}

Back to condition $C(\kappa)$, Theorem~\ref{thm_R*_c*_1} also allows us to refine the first inequality of (\ref{pixroy_hd_nw}) in a natural way.

\begin{corollary}\label{cor_C(kappa)_R*} If $\kappa$ is a cardinal number and $X$ is a topological space with $\omega<\cf(\kappa)<R^*(X)$, then $X$ does not satisfy $C(\kappa)$.

\end{corollary}

\begin{proof} The relations $\cf(\kappa)<R^*(X) =c^*\left(\F[X]\right)$ (see Theorem~\ref{thm_R*_c*_1}) produce a cellular family in $\F[X]$ of cardinality $\cf(\kappa)$. For this reason, Proposition~\ref{prop_weak_precaliber_celullarity} and \ref{prop_caliber_cofinality} imply that $\kappa$ is not a weak precaliber for $\F[X]$. Thus, Theorem~\ref{thm_C(kappa)_equivalence} ensures that $X$ does not satisfy $C(\kappa)$.
\end{proof}

Therefore, the second inequality of (\ref{pixroy_hd_nw}) and Corollary~\ref{cor_C(kappa)_R*} imply that, to determine if $X$ satisfies $C(\kappa)$ or not, we need to focus on those cardinal numbers $\kappa$ with $\cf(\kappa)>\omega$ that satisfy the relations \begin{align}\label{pixroy_R_nw} R^*(X)\leq \cf(\kappa)\leq nw(X).
\end{align}

\begin{question}\label{q_C(kappa)} Is it true that if $\kappa$ is a regular cardinal and $X$ is a $T_1$-space with $R^*(X)\leq\kappa \leq nw(X)$, then $X$ satisfies $C(\kappa)$?

\end{question}

At this point the answer to Question~\ref{q_C(kappa)} depends on the topological space. For example, the ordinal space $[0,\omega_1)$ is right-separated with the natural order; in particular, it is weakly separated. For this reason, Corollary~\ref{cor_thm_R*_c*_2} implies the equalities $$R^*([0,\omega_1))=c\left(\F[[0,\omega_1)]\right)=nw([0,\omega_1))=\omega_1.$$ Thus, $c\left(\F[[0,\omega_1)]\right)=\omega_1$ ensures that $ \F[[0,\omega_1)]$ admits a cellular family of cardinality $\omega_1$ and, therefore, Proposition~\ref{prop_weak_precaliber_celullarity} guarantees that $\F[[0,\omega_1)]$ does not have weak precaliber $\omega_1$; consequently, Theorem~\ref{thm_C(kappa)_equivalence} certifies that $[0,\omega_1)$ {\bf does not} satisfy $C(\omega_1)$.

On the other hand, a classical result states that under Martin's Axiom for $\omega_1$, $\mathsf{MA}(\omega_1)$, any topological space with countable cellularity has precaliber $\omega_1$ (see \cite{roy1989}). Now, since the space $\tilde{\mathbb{R}}$ in Example~\ref{ex_R} satisfies $R^*(\tilde{\mathbb{R}}) = \omega$, Theorem~\ref{thm_R*_c*_2} implies that $\F[\tilde{\mathbb{R}}]$ has countable cellularity and hence, $\mathsf{MA}(\omega_1)$ ensures that $\F[\tilde{\mathbb{R}}]$ has precaliber $\omega_1$. Consequently, Theorem~\ref{thm_C(kappa)_equivalence} indicates that $\tilde{\mathbb{R}}$ {\bf does} satisfy $C(\omega_1)$ under $\mathsf{MA}(\omega_1)$.

\begin{question} Is it possible to prove in \textsf{ZFC} that $\tilde{\mathbb{R}}$ satisfies $C(\omega_1)$? This reduces to giving a proof of the following statement: if $\{x_\alpha : \alpha<\omega_1\} \subseteq \tilde{\mathbb{R}}$ and $\{S_\alpha : \alpha<\omega_1\} \subseteq \mathcal{S}(\mathbb{R})$ satisfy that $S_\alpha$ converges to $x_\alpha$ and $x_\alpha \not \in S_\alpha$ for every $\alpha <\omega_1$, then there exists $J\in [\omega_1]^{\omega_1}$ such that $\{x_\alpha : \alpha\in J\} \cap \bigcup\{S_\alpha : \alpha\in J\} = \emptyset$.

\end{question}

\section{The $C(\kappa)$ condition in topological products}\label{sec_C(kappa)_products}

One question that can be found tacitly in the literature is whether given a pair of spaces $X$ and $Y$ it is satisfied that $\F[X\times Y]$ is homeomorphic to $\F[X]\times \F [Y]$ (see, for example, \cite{vanDouwen1977} and \cite{wage1988}). We will give a way to detect when, for a family of topological spaces $\{X_\alpha : \alpha<\lambda\}$ that fulfills certain properties, the space $\prod_{\alpha<\lambda}\F[X_ \alpha]$ {\bf is not} homeomorphic to $\F\left[\prod_{\alpha<\lambda} X_\alpha\right]$ via its chain conditions (see Theorems~\ref{thm_C(kappa)_product_1}, \ref{thm_C(kappa)_product_2} and \ref{thm_C(kappa)_product_3}).

Regarding the relationship between $C(\kappa)$ and topological products, we start with the following result which follows directly from Lemma~\ref{lemma_C(kappa)_continuous_image} (recall that the natural projections are continuous and surjective).

\begin{proposition}\label{prop_C(kappa)_factors_product} Let $\kappa$ and $\lambda$ be a pair of cardinal numbers ($\lambda$ not necessarily infinite). If $\{X_\alpha : \alpha<\lambda\}$ is a family of topological spaces such that the topological product $\prod\{X_ \alpha : \alpha<\lambda\}$ satisfies $C(\kappa)$, then each factor also satisfies it.
\end{proposition}

Furthermore, the condition $C(\kappa)$ is also preserved under finite products as we will see next.

\begin{proposition}\label{prop_C(kappa)_preservation_two} Let $\kappa$ be a cardinal number. If $X$ and $Y$ are topological spaces that satisfy $C(\kappa)$, then the product $X\times Y$ satisfies it too.

\end{proposition}

\begin{proof} Let $\{W_\alpha : \alpha <\kappa\}$ be a subset of $\tau_{X\times Y}^+$ and for each $\alpha<\kappa$ take $w_\alpha \in W_\alpha$. Let us set collections $\{U_\alpha : \alpha <\kappa\}\subseteq \tau_X$, $\{V_\alpha : \alpha <\kappa\}\subseteq \tau_Y$, $\{x_\alpha : \alpha<\kappa\}\subseteq X$ and $\{y_\alpha : \alpha<\kappa\}\subseteq Y$ such that $U_\alpha\times V_\alpha \subseteq W_\alpha$ and $w_\alpha = (x_\alpha,y_\alpha)$ for each $\alpha<\kappa$. Since $X$ satisfies $C(\kappa)$ there exists $I\in [\kappa]^{\kappa}$ with $\{x_\alpha : \alpha\in I\}\subseteq \bigcap\{ U_\alpha : \alpha \in I\}$. Then, since $Y$ satisfies $C(\kappa)$ there is $J\in [I]^{\kappa}$ such that $\{y_\alpha : \alpha\in J\}\subseteq \bigcap \{V_\alpha : \alpha \in J\}$. Thus, the relation $\{w_\alpha : \alpha \in J\} \subseteq \bigcap\{W_\alpha : \alpha \in J\}$ is fulfilled.
\end{proof}

In this way, an inductive argument can be used to prove the following corollary.

\begin{corollary} Let $\kappa$ be a cardinal number. If $1\leq n<\omega$, $\{X_m : m<n\}$ is a family of topological spaces and each $X_m$ satisfies $C(\kappa)$, then $\prod\{X_m : m<n\}$ satisfies it too.

\end{corollary}

With respect to infinite products, it turns out that it is possible to produce examples to verify that in general the condition $C(\kappa)$ is not preserved. Recall that $t$ denotes the cardinal function known as tightness (see \cite{hodel1984}).

\begin{proposition}\label{prop_C(kappa)_product} Let $\kappa$ and $\lambda$ be a pair of cardinal numbers with $\cf(\kappa)<\lambda$. If $\{X_\alpha : \alpha<\lambda\}$ is a family of Hausdorff spaces with more than one point such that $nw(X_\alpha)<\cf(\kappa)$ for each $\alpha <\lambda$, then each $X_\alpha$ satisfies $C(\kappa)$ and $\prod\{X_\alpha : \alpha<\lambda\}$ does not.

\end{proposition} 

\begin{proof} Let $X:= \prod\{X_\alpha : \alpha<\lambda\}$. On the one hand, since for each $\alpha<\lambda$ we have that $nw(X_\alpha) < \cf(\kappa)$, Proposition~\ref{prop_C(kappa)_nw} implies that $X_\alpha$ satisfies $C(\kappa)$. On the other hand, since $D(2)^{\lambda}$ is embedded in $X$, a routine argument with cardinal functions (see \cite{hodel1984} and \cite{juhasz1980}) shows that $$\lambda = t\left(D(2)^{\lambda}\right) \leq t\left(X\right) \leq hd\left(X\right) \leq nw\left(X\right) = \lambda\cdot \sup\{nw(X_\alpha) : \alpha<\lambda\} = \lambda.$$ Thus, $hd\left(X\right) = \lambda>\cf(\kappa)$ and therefore Corollary~\ref{cor_C(kappa)_hd} implies that $X$ does not satisfy $C(\kappa)$.
\end{proof}

It remains to analyze what is the answer to the question: will it be possible to determine by means of calibers, precalibers and weak precalibers if given a family of spaces $\{X_\alpha : \alpha<\lambda\}$, then the spaces $\prod_{ \alpha<\lambda}\F[X_\alpha]$ and $\F\left[\prod_{\alpha<\lambda} X_\alpha\right]$ are not homeomorphic? We show below that the answer to the previous question is negative for spaces of the form $\F[X\times Y]$ and $\F[X]\times \F[Y]$.

\begin{proposition}\label{prop_C(kappa)_product_two} If $X$ and $Y$ are topological spaces, then the following equalities are true: \begin{align*} \C\left(\F[X\times Y]\right) &= \C\left(\F[X]\right)\cap \C\left(\F[Y]\right); \\
\P\left(\F[X\times Y]\right) &= \P\left(\F[X]\right)\cap \P\left(\F[Y]\right); \ \text{and} \\
\WP\left(\F[X\times Y]\right) &= \WP\left(\F[X]\right)\cap \WP\left(\F[Y]\right).
\end{align*}

\end{proposition}

\begin{proof} For the first equality we note that Theorem~\ref{thm_pixroy_caliber} implies the relations \begin{align*} \C\left(\F[X\times Y]\right) &= \left\{\kappa\in\CN : \cf(\kappa)>|X\times Y|\right\} \\
&= \left\{\kappa\in\CN : \cf(\kappa)>\max\{|X|,|Y|\}\right\} \\
&= \left\{\kappa\in\CN : \cf(\kappa)>|X| \ \wedge \ \cf(\kappa)>|Y|\right\} \\
&= \left\{\kappa\in\CN : \cf(\kappa)>|X| \right\} \cap \left\{\kappa\in\CN : \cf(\kappa)>|Y| \right\} \\
&= \C\left(\F[X]\right)\cap \C\left(\F[Y]\right).
\end{align*}

Now, since the argument for precalibers and weak precalibers is similar, we will only expose the details for precalibers. Naturally, to do this we only have to restrict ourselves to elements of $\UC$. On the one hand, if $\kappa\in \P\left(\F[X\times Y]\right)$ then Theorem~\ref{thm_C(kappa)_equivalence} implies that $X\times Y$ satisfies $C(\kappa)$ and hence Proposition~\ref{prop_C(kappa)_factors_product} ensures that $X$ and $Y$ satisfy $C(\kappa)$. Thus, Theorem~\ref{thm_C(kappa)_equivalence} guarantees that $\F[X]$ and $\F[Y]$ have precaliber $\kappa$. On the other hand, if $\kappa \in \P\left(\F[X]\right)\cap \P\left(\F[Y]\right)$, Theorem~\ref{thm_C(kappa)_equivalence} says that $X$ and $Y$ satisfy $C(\kappa)$. Thus, Proposition~\ref{prop_C(kappa)_preservation_two} certifies that $X\times Y$ satisfies $C(\kappa)$ and therefore Theorem~\ref{thm_C(kappa)_equivalence} asserts that $ \kappa$ is a precaliber for $\F[X\times Y]$.
\end{proof}

Clearly Proposition~\ref{prop_C(kappa)_product_two} can be generalized to any finite product in the natural way.

Now, it is necessary to remember that there are several results in the literature that talk about the preservation of precalibers in topological products. For example, the following theorem compiles some of the work done in \cite{rios2022}, \cite{sanin1948} and \cite{shelah1977}.

\begin{theorem}\label{thm_RSS} Let $\kappa$ be an infinite cardinal, $\lambda$ a cardinal number, $\{X_\alpha : \alpha<\lambda\}$ a family of topological spaces and $X$ a topological space.

\begin{enumerate}
\item If $\kappa$ is regular and uncountable, then $\kappa$ is a caliber (resp., precaliber) for the topological product $\prod\{X_\alpha : \alpha<\lambda\}$ if and only if $\kappa$ is a caliber (resp., precaliber) for each factor.
\item If $\kappa$ is singular, $\cf(\kappa)>\omega$ and $X$ has caliber (resp., precaliber) $\kappa$, then $X^\lambda$ also has caliber (resp., precaliber) $\kappa$.
\end{enumerate}

\end{theorem}

With this result at hand we are ready to prove the following set of theorems.

\begin{theorem}\label{thm_C(kappa)_product_1} Let $\kappa$ and $\lambda$ be a pair of cardinal numbers with $\kappa$ regular and $\omega<\kappa<\lambda$. If $\{X_\alpha : \alpha<\lambda\}$ is a family of Hausdorff spaces with more than one point such that $nw(X_\alpha)<\kappa$ for each $\alpha<\lambda$ , then $$\kappa \in \P\left(\prod_{\alpha<\lambda}\F[X_\alpha]\right) \setminus \WP\left(\F\left[\prod_{\alpha<\lambda } X_\alpha\right]\right).$$

\end{theorem} 

\begin{proof} By Proposition~\ref{prop_C(kappa)_product} every space $X_\alpha$ satisfies $C(\kappa)$ and $\prod\{X_\alpha : \alpha<\lambda\}$ does not. This implies by Theorem~\ref{thm_C(kappa)_equivalence} that each space $\F[X_\alpha]$ has precaliber $\kappa$ and $\F\left[\prod_{\alpha<\lambda } X_\alpha\right]$ does not have weak precalibre $\kappa$. Finally, Theorem~\ref{thm_RSS} guarantees that $\prod_{\alpha<\lambda}\F[X_\alpha]$ does have precaliber $\kappa$.
\end{proof}

The following theorem can be proven analogously to the previous one.

\begin{theorem}\label{thm_C(kappa)_product_2} Let $\kappa$ and $\lambda$ be a pair of cardinal numbers with $\omega<\cf(\kappa)<\lambda$. If $X$ is a Hausdorff space with more than one point such that $nw(X)<\cf(\kappa)$, then  $$\kappa \in \P\left(\F[X]^{\lambda}\right) \setminus \WP\left(\F\left[X^{\lambda}\right]\right).$$

\end{theorem}

\begin{theorem}\label{thm_C(kappa)_product_3} If $X$ is a $T_1$-space and $\lambda$ is a cardinal number (not necessarily infinite), then $$\C\left(\F\left[X^{\lambda}\right]\right) \subseteq \C\left(\F[X]^{\lambda}\right).$$ Furthermore, when $|X|^{\lambda}>|X|$ the previous inclusion is proper.

\end{theorem}

\begin{proof} Let us first note that by Corollary~\ref{prop_weak_precaliber_Hausdorff} and Theorem~\ref{thm_RSS} the relation $\C\left(\F[X]^{\lambda}\right) = \C\left( \F[X]\right) \cap \UC$ holds. Thus, Theorem~\ref{thm_pixroy_caliber} implies that \begin{align*} \C\left(\F\left[X^{\lambda}\right]\right) &= \left\{\kappa \in \CN : \cf(\kappa)>\left|X\right|^{\lambda}\right\} \\
&\subseteq \left\{\kappa \in \CN : \cf(\kappa)>\left|X\right|\right\} = \C\left(\F[X]\right) \cap \UC = \C\left(\F[X]^{\lambda}\right). 
\end{align*}

Finally, if $|X|^{\lambda}>|X|$, then the cardinal $\kappa := |X|^{+}$ satisfies the relations $\kappa\in \C\left(\F[ X]^{\lambda}\right)$ and $\kappa \not \in \C\left(\F\left[X^{\lambda}\right]\right)$.
\end{proof}

In particular, Theorems ~\ref{thm_C(kappa)_product_1}, \ref{thm_C(kappa)_product_2}, and \ref{thm_C(kappa)_product_3} show that chain conditions can also be used to detect when $ \prod_{\alpha<\lambda}\F[X_\alpha]$ and $\F\left[\prod_{\alpha<\lambda} X_\alpha\right]$ are not homeomorphic.

\section{The weak Lindel\"of degree of $\F[X]$}

If $X$ is a topological space, then the {\it weak Lindel\"of degree} of $X$, $wL(X)$, is the cardinal number $$\min\left\{\kappa\in\CN : \forall \ \mathcal{U}\subseteq\tau_X\left(X=\bigcup\mathcal{U}\to \exists\mathcal{V}\in[\mathcal{U}]^{\leq\kappa}\left(X=\overline{\bigcup\mathcal{V}}\right)\right)\right\}.$$ We shall say that $X$ is {\it weakly Lindel\"of} if $wL(X)=\omega$. Additionally, we will work with the cardinal functions $$L^*(X) := \sup \{L(X^{n}) : n\in \mathbb{N}\} \quad \text{and} \quad hd_c(X) := \sup\{d(A) : A\in \tau_X^{*}\}.$$ Clearly, $d(X) \leq hd_c(X) \leq hd(X)$.

In \cite[Question~3.23, p.~3087]{sakai2012} Sakai asked whether the fact that $\F[X]$ is weakly Lindel\"of implies that $X$ is hereditarily separable and proved that if $X $ is countably tight then the previous question has an affirmative answer. We shall expand Sakai's result by proving that if $\F[X]$ is weakly Lindel\"of and $X$ satisfies any of the following conditions: \begin{itemize}
 \item $X$ is a Hausdorff $k$-space;
 \item $X$ is a countably tight $T_1$-space;
 \item $X$ is weakly separated,
 \end{itemize} then $X$ is hereditarily separable (see Corollary~\ref{cor_wL_k-space_2}).

What follows is intended to show that if $X$ is a countably tight $T_1$-space or $X$ is a Hausdorff $k$-space, then \begin{align}\label{wL_relation} L^*(X)\cdot hL(X) \cdot hd(X) \leq wL\left(\F[X]\right).
\end{align} We will prove (\ref{wL_relation}) by establishing generalizations and connections between various results of Sakai and Tall exposed in \cite{sakai2012} and \cite{tall1974} respectively.

If $\mathcal{C}$ is a collection of subsets of $X$, we say that $\mathcal{C}$ is an {\it $\omega$-cover} for $X$ if for any $F\in [X]^{<\omega}$ there exists $C\in \mathcal{C}$ with $F\subseteq C$. In \cite{gernag1982} Gerlits and Nagy showed that if $X$ is a topological space, then $X^n$ is Lindel\"of for every $n\in \mathbb{N}$ if and only if all $\omega$-open covers of $X$ admit a countable $\omega$-subcover. The following lemma generalizes this fact and its proof can be found in \cite[S.~148, p.~122]{tkachuk2011}.

\begin{lemma}\label{lemma_wL_L*} If $X$ is a topological space and $\kappa$ is an infinite cardinal, then $L^*(X) \leq \kappa$ if and only if every open $\omega$-cover of $X$ admits an $\omega$-subcover of size at most $\kappa$.

\end{lemma}

A family $\mathcal{A}$ formed by subsets of a space $X$ is {\it discrete} if for any $x\in X$ there exists $U\in \tau_X$ with $x\in U$ and $ |\{A\in \mathcal{A} : A\cap U \neq \emptyset\}|\leq 1$ (see \cite[p.~16]{engelking1989}). The {\it discrete cellularity} of $X$ will be the cardinal number $$dc(X) := \sup\left\{|\mathcal{U}| : \mathcal{U}\subseteq \tau_X^+ \ \text{is a discrete family}\right\} +\omega.$$ According to \cite[Definition~3.1, p.~3083]{sakai2012}, we say that $X$ satisfies the {\it discrete countable chain condition} (in symbols, dccc) if $dc(X)=\omega$.

An immediate observation is that if $X$ is a topological space and $A$ is a clopen subset of $X$, then $dc(A) \leq dc(X)$. Also, any discrete family $\mathcal{U}\subseteq \tau_X^+$ is a cellular family; consequently, it is satisfied that $dc(X)\leq c(X)$. Furthermore, simple reasoning shows that $wL(X) \leq c(X)$ (see \cite[p.~16]{hodel1984}). In the realm of $T_3$-spaces we can also connect $dc(X)$ with $wL(X)$.

\begin{lemma}\label{lemma_wL_dc} If $X$ is a $T_3$-space, then $dc(X) \leq wL(X)$.

\end{lemma}

\begin{proof} Set $\kappa := wL(X)$. Let $\lambda$ be an infinite cardinal, $\mathcal{U}:=\{U_\alpha : \alpha<\lambda\} \subseteq \tau_X^+$ a discrete family, and for each $\alpha<\lambda$ let $x_\alpha \in U_\alpha$ and $V_\alpha \in \tau_X$ satisfy $x_\alpha \in V_\alpha \subseteq \overline{V_\alpha } \subseteq U_\alpha$. Observe that, since $\{U_\alpha : \alpha<\lambda\}$ is a discrete family, then $\{\overline{V_\alpha} : \alpha<\lambda\}$ is also a discrete family; in particular, since $\{\overline{V_\alpha} : \alpha<\lambda\}$ is locally finite it satisfies that $V:= X\setminus \bigcup_{\alpha<\lambda} \overline{V_\alpha}$ is an open subset of $X$ (see \cite[Corollary~1.1.12, p.~17]{engelking1989}).

Thus, the collection $\mathcal{U}\cup \{V\}$ is an open cover of $X$ and thus there exists $\mathcal{W} \in [\mathcal{U}\cup \{V \}]^{\leq \kappa}$ with $\overline{\bigcup\mathcal{W}}=X$. Our goal now is to check that $\mathcal{U}$ is a subset of $\mathcal{W}$. It turns out that if $\alpha<\lambda$ then there exists $W\in \mathcal{W}$ with $V_\alpha \cap W \neq \emptyset$. Then, since $V_\alpha$ and $V$ are disjoint sets, necessarily $W\in \mathcal{U}$. Finally, since $V_\alpha \subseteq U_\alpha$ and $\mathcal{U}$ is a cellular family, it is verified that $U_\alpha = W$ and, therefore, that $U_\alpha \in \mathcal{W}$. In conclusion, $|\mathcal{U}|\leq \kappa$.
\end{proof}

The following pair of lemmas can be proven with routine arguments (the second one can be found in \cite[Lemma~3.9, p.~3084]{sakai2012}).

\begin{lemma}\label{lemma_wL_hL_L(U)} If $X$ is a topological space, then $hL(X) = \sup\{L(U) : U\in \tau_X\}$.

\end{lemma}

\begin{lemma}\label{lemma_wL_clopen} If $X$ is a topological space and $\mathcal{U}$ is a subset of $\tau_X$, then $\mathcal{V}(\mathcal{U}) :=\{F\in \F[X ] : \exists U\in\mathcal{U}(F\subseteq U)\}$ is a clopen subset of $\F[X]$.

\end{lemma}

What follows is to generalize the first part of \cite[Theorem~3.10, p.~3085]{sakai2012}.

\begin{theorem}\label{thm_wL_L*_hL_dc} If $X$ is a topological space, then $\max\{L^*(X), hL(X)\}\leq dc\left(\F[X]\right)$.

\end{theorem}

\begin{proof} Set $\kappa := dc\left(\F[X]\right)$. To see that $L^*(X) \leq \kappa$, let $\{U_\alpha : \alpha<\lambda\}$ be an open $\omega$-cover of the space $X$. For each $\alpha<\lambda$ consider the set $$V_\alpha := \mathcal{V}\left(\{U_\alpha\})\setminus \mathcal{V}(\{U_\beta : \beta<\alpha\}\right).$$ By Lemma~\ref{lemma_wL_clopen} it is satisfied that $\mathcal{V}:=\{V_\alpha : \alpha<\lambda\}$ is a family of clopen subsets of $\F[X]$. Furthermore, if $F\in \F[X]$ and $\beta := \min \{\alpha<\lambda : F\subseteq U_\alpha\}$, then $F\in V_\beta$; that is, $\mathcal{V}$ is a cover of $\F[X]$. Furthermore, if $\alpha<\beta<\lambda$ and we assume that $F\in V_\alpha \cap V_\beta$ exists, then necessarily $F \subseteq U_\alpha$ and $F\setminus U_\alpha \neq \emptyset$, a contradiction. Therefore, $\mathcal{V}$ is a pairwise disjoint cover of $\F[X]$ made up of open and closed subsets of $\F[X]$.

Now, if $J:= \{\alpha<\lambda : V_\alpha \neq \emptyset\}$, then since $\{V_\alpha : \alpha\in J\}$ is a discrete family in $\ F[X]$ it is true that $|J| \leq \kappa$. Thus, $\{U_\alpha : \alpha\in J\}$ is an $\omega$-subcover of $\{U_\alpha : \alpha<\lambda\}$ of size at most $\kappa$. Indeed, if $F\in \F[X]$ then there exists $\alpha \in J$ with $F \in V_\alpha$ and therefore $F \subseteq U_\alpha$. By virtue of Lemma~\ref{lemma_wL_L*}, this argument proves that $L^*(X) \leq \kappa$.

To verify that $hd(X) \leq \kappa$ we observe that if $U$ is an open subset of $X$, then the identity function is a homeomorphism between $\F[U]$ and $\mathcal{V} (\{U\})$. Thus, since $\mathcal{V}(\{U\})$ is an open and closed subset of $\F[X]$ it satisfies that $dc\left(\F[U]\right) \leq \kappa$ and so, the first part of this result guarantees that $L^*(U) \leq \kappa$; in particular, $L(U)\leq \kappa$. In conclusion, Lemma~\ref{lemma_wL_hL_L(U)} ensures that $hL(X) \leq \kappa$.
\end{proof}

Since $\F[X]$ is a zero-dimensional Hausdorff space when $X$ is $T_1$, Lemma~\ref{lemma_wL_dc} and Theorem~\ref{thm_wL_L*_hL_dc} produce the following result.

\begin{theorem}\label{thm_wL_L*_hL} If $X$ is a $T_1$-space, then $L^*(X)\cdot hL(X) \leq wL\left(\F[X]\right)$.

\end{theorem}

We now trace the path to prove that $hd(X) \leq wL\left(\F[X]\right)$ if $X$ is a Hausdorff $k$-space. What follows is to expose some generalizations of the results of the third section of \cite{tall1974}.

\begin{proposition}\label{prop_cf_density_tall} If $X$ is a topological space and $d(X)=|X|$, or $t(X)<\cf\left(d(X)\right)$, then $\cf\left( d(X)\right)$ and $d(X)$ are not calibers for $X$.

\end{proposition}

\begin{proof} Since both arguments differ slightly, we can do both cases simultaneously. Let's put $\kappa:=d(X)$ and $\lambda:=\cf\left(d(X)\right)$. Let $\{x_\beta: \beta<\kappa\}$ be a subset of $X$ such that $X=\overline{\{x_\beta : \beta<\kappa\}}$ (in the first case, we will also assume that $X=\{x_\beta : \beta<\kappa\}$). Finally, let $f:\lambda \to \kappa$ be a strictly increasing and cofinal function.

\medskip

\noindent {\bf Claim.} For every $J\in [\lambda]^{\lambda}$ it is satisfied that $X= \bigcup_{\alpha\in J}\overline{\left\{x_\beta : \beta < f(\alpha)\right\}}$.

\medskip

Fix $J\in [\lambda]^{\lambda}$ and $x\in X$. Under these circumstances  $\restr{f}{J}:J\to\kappa$ is also cofinal in $\kappa$. Now, when $d(X)=|X|$, there exists $\gamma<\kappa$ such that $x=x_\gamma$; thus, if $\alpha\in J$ is such that $\gamma<f(\alpha)$, then $x\in \{x_\beta : \beta < f(\alpha)\}$ holds. On the other hand, if $t(X)<\cf\left(d(X)\right)$, we use $x\in \overline{\{x_\beta : \beta<\kappa\}}$ for find $I\in [\kappa]^{\leq t(X)}$ such that $x\in \overline{\{x_\beta : \beta\in I\}}$. Then, since $I$ is not cofinal in $\kappa$, there exists $\alpha<\lambda$ such that $I\subseteq [0,f(\alpha))$ and thus, we deduce that $x\in \overline {\{x_\beta : \beta < f(\alpha)\}}$.

Now, for each $\alpha<\lambda$ define $U_\alpha := X\setminus \overline{\left\{x_\beta : \beta < f(\alpha)\right\}}$. It turns out that $\{U_\alpha : \alpha<\lambda\}$ is a subset of $\tau_X^+$. By virtue of the Claim above, the set $\bigcap\{U_\alpha : \alpha\in J\}$ is empty when $J\in [\lambda]^{\lambda}$. Consequently, $X$ does not have caliber $\lambda$ and thus, by Proposition~\ref{prop_caliber_cofinality}, $\kappa$ is also not a caliber for $X$.
\end{proof}

A cardinal function $\phi$ {\it reflects} a cardinal number $\kappa$ if the condition $\phi(X)\geq\kappa$ implies the existence of $Y\in[X]^{\leq\kappa}$ with $\phi(Y)\geq \kappa$. For example, a classical result states that density reflects any regular cardinal, that is, if $\kappa$ is a regular cardinal and $d(X)\geq \kappa$, then there exists $Y\subseteq X$ with $ |Y|=d(Y)= \kappa$ (see \cite[Theorem~2.5, p.~54]{hodvau2000}). The next lemma follows from this fact.

\begin{lemma}\label{lemma_reflection_tall} If $X$ is a topological space such that $hd(X)>\kappa$, then there exists $Y\subseteq X$ with $|Y|=d(Y)= \kappa^{+}$.

\end{lemma}

It is not difficult to check that if $D$ is a dense subspace of $X$, then $d(X)\leq d(D)\leq d(X)\cdot t(X)$ (see \cite[Theorem~ 3.9, p.~17]{hodel1984}). In particular, for any subset $Y$ of $X$ it holds that \begin{align}\label{relation_density_tightness} d(\overline{Y})\leq d(Y)\leq d(\overline{Y})\cdot t(\overline{Y}).
\end{align} We will use these relationships in our next result that generalizes to \cite[Theorem~3.24, p.~324]{tall1974}.

\begin{theorem}\label{thm_tall} Let $X$ be a topological space and $\kappa$ a cardinal number such that $t(X)\leq \kappa$. If all closed subspaces of $X$ have caliber $\kappa^{+}$, then $hd(X)\leq\kappa$.

\end{theorem}

\begin{proof} If $hd(X)>\kappa$, then Lemma~\ref{lemma_reflection_tall} produces $Y\subseteq X$ such that $d(Y)= \kappa^{+}$. Since the property {\lq\lq}$t(\boldsymbol\cdot)\leq \kappa${\rq\rq} is hereditary, the inequalities of (\ref{relation_density_tightness}) guarantee that $d(\overline{Y })\leq \kappa^{+}\leq d(\overline{Y})\cdot \kappa$ and, therefore, we obtain the equality $d(\overline{Y})=\kappa^{+}$ . Thus, since $t(\overline{Y})<\cf\left(d(\overline{Y})\right)$, Proposition~\ref{prop_cf_density_tall} implies that $\overline{Y}$ does not have caliber $\kappa^{+}$.
\end{proof}

Recall that a topological space $X$ is a {\it $k$-space} if for any $A\subseteq X$ it is satisfied that $A$ is closed in $X$ provided that for any compact subset $K$ of $X$, $A\cap K$ is closed in $K$. In \cite{sapirovski1972} \v{S}apirovski\u{\i} proved that if $X$ is a Hausdorff $k$-space, then $t(X)\leq s(X)$. Our next corollary follows from this result.

\begin{corollary}\label{cor_thm_tall} If $X$ is a Hausdorff $k$-space, $\kappa$ is a cardinal number, and all closed subspaces of $X$ have caliber $\kappa^+$, then $hd(X)\leq \kappa $.

\end{corollary}

\begin{proof} Under our hypotheses it is satisfied that all the closed subspaces of $X$ have cellularity less than $\kappa^+$; consequently, all its subspaces have cellularity less than $\kappa^+$ and, therefore, $s(X)\leq \kappa$. Finally, since $t(X)\leq \kappa$, Theorem~\ref{thm_tall} guarantees that $hd(X)\leq \kappa$.
\end{proof}

The following result is a generalization of the first part of \cite[Proposition~3.13, p.~3085]{sakai2012}.

\begin{lemma}\label{lemma_wL_closed} If $X$ is a topological space and $A$ is a closed subspace of $X$, then $d(A) \leq wL\left(\F[X]\right)$. Consequently, $hd_c(X) \leq wL\left(\F[X]\right)$.

\end{lemma}

\begin{proof} Set $\kappa:=wL\left(\F[X]\right)$. Since $\{[\{x\},X] : x\in A\} \cup \{[F,X\setminus A] : F\in \F[X] \ \text{y} \ F\cap A = \emptyset\}$ is an open cover for $\F[X]$, there are $\{x_\alpha : \alpha<\kappa\} \subseteq A$ and $\{F_\alpha : \alpha <\kappa\} \subseteq \F[X]$ with $A \cap \bigcup\{F_\alpha : \alpha<\kappa\}= \emptyset$ and such that $\mathcal{U} := \{[\{x_\alpha\},X] : \alpha<\kappa\} \cup \{[F_\alpha,X\setminus A] : \alpha<\kappa\}$ satisfies $\F[X ] = \overline{\bigcup\mathcal{U}}$. To see that $D:= \{x_\alpha : \alpha<\kappa\}$ is dense in $A$ suppose, in search of a contradiction, that there exists $x\in A \setminus \overline{D}$ . It turns out that if $U := A \setminus \overline{D}$, then $[\{x\},U]$ is a non-empty open subset of $\F[X]$ that has empty intersection with every $[ F_\alpha,X\setminus A]$; consequently, since $\F[X] = \overline{\bigcup\mathcal{U}}$, there exists $\alpha<\kappa$ with $[\{x\},U] \cap [\{x_\alpha\},X] \neq \emptyset$. However, the latter implies that $x_\alpha$ is not an element of $\overline{D}$, a contradiction. In conclusion, $d(A) \leq \kappa$.
\end{proof}

\begin{theorem}\label{thm_wL_L*_hL_wL} If $X$ is a Hausdorff $k$-space or $X$ has countable tightness, then $hd(X)\leq wL\left(\F[X]\right)$.

\end{theorem}

\begin{proof} Set $\kappa:=wL\left(\F[X]\right)$. If $X$ has countable tightness, then from (\ref{relation_density_tightness}) it follows that $hd(X) = hd_c(X)$. Thus, Lemma~\ref{lemma_wL_closed} ensures that $hd(X)\leq \kappa$. Now, if $X$ is a Hausdorff $k$-space, then Lemma~\ref{lemma_wL_closed} implies that any closed subspace of $X$ has density at most $\kappa$; consequently, Proposition~\ref{prop_caliber_density} guarantees that all closed subspaces of $X$ have caliber $\kappa^+$ and hence, Corollary~\ref{cor_thm_tall} certifies that $hd(X)\leq \kappa$.
\end{proof}

In conclusion, Theorems~\ref{thm_wL_L*_hL} and \ref{thm_wL_L*_hL_wL} produce the following corollaries.

\begin{corollary}\label{cor_wL_k-space_1} If $X$ is a Hausdorff $k$-space or $X$ is a $T_1$-space with countable tightness, then $$L^*(X)\cdot hL(X) \cdot hd(X) \leq wL\left(\F[X]\right).$$

\end{corollary}

The connections between some of the cardinal functions of $X$ and $\F[X]$ that we have established so far are summarized in the following diagram (the relation $\kappa\to \lambda$ means $\kappa\geq \lambda$):

\begin{figure}[H]
\begin{center}
\begin{tikzcd}

c\left(\F[X]\right) \arrow[rr, bend left=30]{}{} \arrow[dd]{}{} & hd_c(X) & wL\left(\F[X]\right) \arrow[dd]{}{} \arrow[ddl]{}{} \arrow[dl,swap]{}{(*)} \arrow[l]{}{} \arrow[r]{}{} & dc\left(\F[X]\right)\\

& hd(X) \arrow[u]{}{} & & \\

R(X) \arrow[r]{}{} \arrow[ru]{}{}& hL(X) & L^*(X) & \\

\end{tikzcd}
\end{center}
\vspace{-10mm}
\caption{Relations between some cardinal functions of $X$ and $\F[X]$}\label{figura_pixroy}
\end{figure}
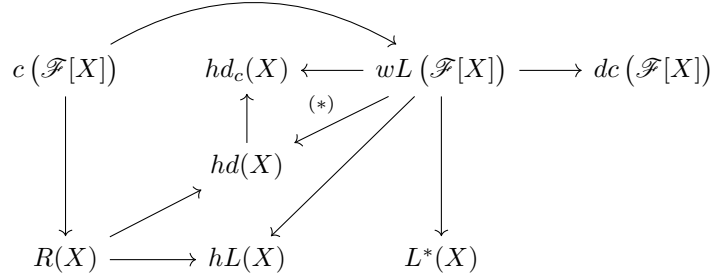
Naturally, in the relation $wL\left(\F[X]\right)\to hd(X)$ we need $X$ to be $T_1$ and have countable tightness, or else to be a Hausdorff $k$-space.

In general, we don't know if it is possible to directly connect $R(X)$ and $dc\left(\F[X]\right)$ in some way, but we can establish a relationship in this regard by considering a variant for the $R$ function. The following result will be essential to achieve this objective.

\begin{lemma}\label{lemma_wL_weak_sep} If $X$ is weakly separated, then $R(X) = dc\left(\F[X]\right)=|X|$.

\end{lemma}

\begin{proof} Clearly, $R(X) = |X| \geq dc\left(\F[X]\right)$. To see that $dc\left(\F[X]\right) \geq |X|$, let $\{U_x : x\in X\}$ be a weak separation for $X$ and let us check that $\{[\{x\},U_x] : x \in X\}$ is a discrete family in $\F[X]$. For $F\in \F[X]$ define $U_F:= \bigcup\{U_x : x\in F\}$ and suppose, in search of a contradiction, that there are distinct $x,y\in X$ such that $[F,U_F] \cap [\{x\},U_x] \neq \emptyset \neq [F,U_F] \cap [\{y\},U_y]$. It then turns out that $F \subseteq U_x\cap U_y$ and $z,w\in F$ exist with $x\in U_z$ and $y\in U_w$. Now, if $x\neq z$ or $y\neq w$, then $x\in U_z$ and $z\in U_x$, or $y\in U_w$ and $w\in U_y$; which is absurd. Thus, $x=z$ and $y=w$; in particular, $\{x,y\}\subseteq F$ and hence $x\in U_y$ and $y\in U_x$, a contradiction. This argument shows that $\{[\{x\},U_x] : x \in X\}$ is a discrete family in $\F[X]$ of cardinality $|X|$.
\end{proof}

If $X$ is a weakly separated $T_1$-space, Lemmas~\ref{lemma_wL_dc} and~\ref{lemma_wL_weak_sep} imply that $R(X) = wL\left(\F[X]\right)= |X|$. However, it is also possible to obtain the same result without the additional assumption of being $T_1$.

\begin{lemma}\label{lemma_wL_weak_sep_2} If $X$ is weakly separated, then $R(X) = wL\left(\F[X]\right)=|X|$.

\end{lemma}

\begin{proof} It is enough to show that $wL\left(\F[X]\right) \geq |X|$. Let $\{U_x : x\in X\}$ be a weak separation for $X$ and for each $F\in \F[X]$ define $U_F:= \bigcup\{U_x : x\in F\} $. Clearly, $\{[F,U_F] : F\in \F[X]\}$ is an open cover of $X$. Let $\mathcal{F}$ be a subset of $\F[X]$ with $|\mathcal{F}|<|X|$ and $x\in X \setminus \bigcup\mathcal{F}$. It turns out that for each $F\in \mathcal{F}$ the relation $[x,U_x] \cap [F,U_F] = \emptyset$ is satisfied. Indeed, if this intersection were non-empty, then $x\in U_F$ and so there would exist $y\in F$ with $x\in U_y$. On the other hand, as $F \subseteq U_x$ it also happens that $y\in U_x$, a contradiction to weak separation. This argument proves that $\{[F,U_F] : F\in \F[X]\}$ does not admit subfamilies of cardinality less than $|X|$ with dense union.
\end{proof}

\begin{corollary}\label{cor_wL_k-space_2} Consider the following conditions for a topological space $X$.

\begin{enumerate}
\item $X$ is a Hausdorff $k$-space.
\item $X$ has countable tightness and is $T_1$.
\item $X$ is weakly separated.
\end{enumerate}

If $X$ satisfies any of the conditions from (1) to (3) and $\F[X]$ is weakly Lindel\"of, then $X$ is hereditarily separable.

\end{corollary}

\begin{proof} If $X$ satisfies (1) or (2) the result is a consequence of Corollary~\ref{cor_wL_k-space_1}. Finally, if $X$ is weakly separated, since $\F[X]$ is weakly Lindel\"of Lemma~\ref{lemma_wL_weak_sep_2} guarantees that $X$ is even countable.
\end{proof}

Now, consider the cardinal function $R_o(X)$ determined by the rule $$\sup\{|U| : U \ \text{is a weakly separated open subspace of} \ X\}+\omega.$$ Naturally, $R_o(X) \leq R(X)$. On the other hand, recall that if $U\in \tau_X$ then $\F[U]$ is a clopen subspace of $\F[X]$. We will use this fact in the proof of the last result of this section.

\begin{theorem} If $X$ is a topological space, then $R_o(X) \leq dc\left(\F[X]\right)$.

\end{theorem}

\begin{proof} Let $U$ be a weakly separated open subspace of $X$. Observe that by Lemma~\ref{lemma_wL_weak_sep} it is satisfied that $dc\left(\F[U]\right) = |U|$. Also, since $\F[U]$ is a clopen subspace of $\F[X]$, it follows that $dc\left(\F[U]\right) \leq dc\left(\F[X] \right)$. Thus, $|U| \leq dc\left(\F[X]\right)$ and consequently $R_o(X) \leq dc\left(\F[X]\right)$. 
\end{proof}

One of the natural questions that arises is whether the functions $R(X)$ and $R_o(X)$ are different. Fortunately, one of the reviewers suggested the following example. For an infinite cardinal $\kappa$, let $D(2)^{\kappa}$ stand for the Cantor cube of weight $\kappa$. On the one hand, since the discrete space $D(\kappa)$ embeds into $D(2)^{\kappa}$ and $nw\left(D(2)^{\kappa}\right)=\kappa$ (see \cite{hodel1984}), it follows that $R(D(2)^{\kappa})=\kappa$ (see Proposition~\ref{prop_weak_sep_car_fun}). On the other hand, since the Cantor space $D(2)^{\omega}$ is an uncountable cosmic space that embeds into every non-empty open subset of $D(2)^{\kappa}$, it is derived that $R_o(D(2)^{\kappa})=\omega$ (see Corollary~\ref{cor_weak_sep_uncountable}).

However, the following questions remain to be answered:

\begin{question} Let $X$ be a topological space.

\begin{enumerate}
\item What is the relationship between $R(X)$ and $dc\left(\F[X]\right)$?
\item Is it possible to connect $R(X)$ and $L^*(X)$?
\end{enumerate}

\end{question}

The author thanks Dr. \'Angel Tamar\'iz Mascar\'ua for the constant exchange of ideas and guidance in the preparation of this work. In addition, he also thanks the reviewers for their thoughtful comments regarding this paper.


\begin{thebibliography}{99}
  \thispagestyle{myheadings}
  
  \bibitem{creede1970} G. D. Creede, \textit{Concerning semistratifiable spaces}, Pacific J. Math., {\bf 32} (1970), 47--54.    
  
  \bibitem{vanDouwen1977} E. van Douwen, \textit{The Pixley-Roy topology on spaces of subsets}, Set-theoretic topology (G. M. Reed, ed.), pp. 111--134, 1977.    
  
  \bibitem{engelking1989} R. Engelking, \textit{General Topology}, Sigma series in pure mathematics, Heldermann, vol. 6, 1989.
  
  \bibitem{gernag1982} J. Gerlits, Z. S. Nagy, \textit{Some properties of $C(X), I$}, Topol. Appl., {\bf 14(2)} (1982), 151--161. 
  
  \bibitem{hodel1984} R. Hodel, \textit{Cardinal functions I}, Handbook of set-theoretic topology (K. Kunen y J. E. Vaughan, eds.), pp. 1--61, 1984. 
  
    \bibitem{hodvau2000} R. E. Hodel, J. E. Vaughan, \textit{Reflection theorems for cardinal functions}, Topol. Appl., {\bf 100} (2000), 47--66. 
  
  \bibitem{juhasz1980} I. Juh\'asz, \textit{Cardinal Functions in Topology - Ten Years Later}, Math. Centre Tracts, vol. 123, 1980.
  
  \bibitem{juhshe2003} I. Juh\'asz, S. Shelah, \textit{Generic left-separated spaces and calibers}, Topol. Appl., {\bf 132} (2003), 103--108.
  
  \bibitem{kunen1980} K. Kunen, \textit{Set Theory. An introduction to independence proofs}, Studies in logic and the fundations of mathematics, North-Holland publishing co., vol. 102, 1980. 
  
  \bibitem{lutzer1978} D. J. Lutzer, \textit{Pixley-Roy topology}, Topol. Proc., {\bf 3} (1978), 139--158.
  
  \bibitem{mcintyre2006} D. W. McIntyre, {\it A regular countable chain condition space without compact-caliber $(\omega_1,\omega)$}, Annals of the New York Academy of Sciences {\bf 704(1)} (2006), 269--272.
 
 \bibitem{miller2017} A. W. Miller, \textit{Descriptive Set Theory and Forcing: How to Prove Theorems about Borel Sets the Hard Way}, Lecture Notes in Logic, Cambridge University Press, 2017.
  
  \bibitem{pixroy1969} C. Pixley, P. Roy, \textit{Uncompletable Moore spaces}, Proc. Auburn Univ. Conf. (W. R. R. Transue, ed.), pp. 75--85, 1969.
  
  \bibitem{rios2022} A. R\'ios-Herrej\'on, \textit{Singular precalibers for topological products}, Topol. Appl., {\bf 317} (2022), 108190.
  
 \bibitem{riotam2023} A. R\'ios-Herrej\'on, \'A. Tamariz-Mascar\'ua, \textit{Some notes on topological calibers}, Preprint, \url{http://arxiv.org/abs/2302.12408}, 2023.

  \bibitem{roy1989} N. M. Roy, \textit{Is the product of ccc spaces a ccc space?}, Publicacions Matem\`atiques, {\bf 33} (1989), 173--183.
  
  \bibitem{sakai2012} M. Sakai, \textit{Cardinal functions of Pixley-Roy hyperspaces}, Topol. Appl., {\bf 159} (2012), 3080--3088.
  
  \bibitem{sakai1983} S. Sakai, \textit{Cardinal functions on Pixley-Roy hyperspaces}, Proc. Amer. Math. Soc., {\bf 89(2)} (1983), 336--340. 
  
  \bibitem{sanin1948} N. A. \v{S}anin, {\it On the product of topological spaces}, Trudy Mat. Inst. Steklov., {\bf 24} (1948), (In Russian), 1--112.
  
  \bibitem{sapirovski1972} B. \v{S}apirovski\u{\i}, \textit{On discrete subspaces of topological spaces; weight, tightness and Suslin number}, Soviet Math. Dokl., {\bf 13} (1972), 215--219.
  
 \bibitem{shelah1977} S. Shelah, \textit{Remarks on cardinal invariants in topology}, Gen. Topol. Appl., {\bf 7(3)} (1977), 251--259.
 
   \bibitem{tall1974} F. D. Tall, \textit{The countable chain condition versus separability --- Applications of Martin's axiom}, Gen. Topol. Appl., {\bf 4 }(1974), 315--339.
  
  \bibitem{tkachenko1978} M. G. Tkachenko, \textit{Chains and cardinals}, Dokl. Akad. Nauk SSSR, {\bf 239(3)} (1978), Soviet Math. Dokl., {\bf 19(2)} (1978) (English translation), 382--385.
  
  \bibitem{tkachuk2011} V. V. Tkachuk, \textit{A $C_p$-Theory Problem Book. Topological and function spaces}, Problem books in mathematics, Springer-Verlag New York, 2011. 
  
  \bibitem{wage1988} M. L. Wage, \textit{Homogeneity of Pixley-Roy spaces}, Topol. Appl, {\bf 28} (1988), 45--57.
   
\end{thebibliography}
\end{document}